\documentclass[12pt]{article}

\usepackage{amsthm,color,amsmath,amsfonts,amssymb}

\parskip=\medskipamount
\def\onehalf{{\textstyle\frac12}}

\def\cinfty#1{C^{\scriptscriptstyle\infty}(#1)}
\def\vectorfields#1{{\mathcal X}(#1)}
\def\oneforms#1{{\mathcal X}^*(#1)}

\def\tvectorfields{\vectorfields{\tau}}

\def\del{{\nabla}}
\def\lie#1{{\mathcal L}_{#1}}
\def\fpd#1#2{\frac{\partial #1}{\partial #2}}
\def\R{{\rm I\kern-.20em R}}
\def\sode{{\sc sode}}

\newtheorem{prop}{\bf Proposition}
\newtheorem{lemma}{\bf Lemma}

\newtheorem{dfn}{\bf Definition}

\def\DV#1{\V{\rm D}_{#1}}
\def\DH#1{\H{\rm D}_{#1}}

\def\H#1{{#1}^{\scriptscriptstyle H}}
\def\V#1{{#1}^{\scriptscriptstyle V}}

\def\tr{\mathop{\rm tr}}
\def\cof{\mathop{\rm cof}}
\def\jb{\bar{J}_2{}}

\def\hook{{\mathchoice
{\vrule height 0pt depth 0.4pt width 3pt \vrule height 5pt depth
0.4pt \kern 3pt} {\vrule height 0pt depth 0.4pt width 3pt  \vrule
height 5pt depth 0.4pt \kern 3pt} {\vrule height 0pt depth 0.2pt
width 1.5pt  \vrule height 3pt depth 0.2pt width 0.2pt \kern 1pt}
{\vrule height 0pt depth 0.2pt width 1.5pt  \vrule height 3pt depth
0.2pt width 0.2pt \kern 1pt} }}

\begin{document}

\title{Driven cofactor systems and Hamilton-Jacobi separability}

\author{W.\ Sarlet$^{a,b}$ and G.\ Waeyaert$^{a}$\\
{\small ${}^a$Department of Mathematics, Ghent University }\\
{\small Krijgslaan 281, B-9000 Ghent, Belgium}\\[1mm]
{\small ${}^b$Department of Mathematics and Statistics, La Trobe University}\\
{\small Bundoora, Victoria 3086, Australia}
}
\date{}
\maketitle

\begin{quote}
{\bf Abstract.} {\small This is a continuation of the work initiated
in \cite{SarBier} on so-called driven cofactor systems, which are
partially decoupling second-order differential equations of a
special kind. The main purpose in \cite{SarBier} was to obtain an
intrinsic, geometrical characterization of such systems, and to
explain the basic underlying concepts in a brief note. In the
present paper we address the more intricate part of the theory. It
involves in the first place understanding all details of an
algorithmic construction of quadratic integrals and their
involutivity. It secondly requires explaining the subtle way in
which suitably constructed canonical transformations reduce the
Hamilton-Jacobi problem of the (a priori time-dependent) driven part
of the system into that of an equivalent autonomous system of
St\"ackel type. }
\end{quote}

\section{Introduction}

To the best of our knowledge, the idea of a cofactor system stems
from a paper by Rauch-Wojciechowski {\it et al.\/} on certain
Newtonian systems in the Euclidean plane \cite{WML}, but the term
`quadratic integral of cofactor type' came in use with Lundmark's
generalization to systems of arbitrary dimension, which was
published with some delay in \cite{Lundmark2}. The term refers to
the fact that the matrix of the quadratic part of the first integral
comes from the cofactor tensor of a tensor which has special
properties with respect to the (Euclidean) metric of the kinetic
energy of the system (a tensor which was termed `inertia tensor' in
Benenti's work \cite{Benenti} on St\"ackel systems). The other point
to be emphasized about the `quadratic integral of cofactor type' is
that the zeroth-order terms in the integral do not come from a
potential energy function, in other words, the Newtonian systems
under consideration have force terms which are of nonconservative
type, albeit of a very special nature, determined by a scalar
function also and by the cofactor tensor. It was recognized in
\cite{CS2001} that the work of Lundmark could easily be generalized
to kinetic energy terms associated to an arbitrary Riemannian metric
and in fact fits perfectly within the theory of `special conformal
Killing tensors', as studied in \cite{CST}. This gave rise to an
intrinsic characterization of general cofactor systems, in the
context of which it was shown that the corresponding quadratic
integral in fact is the Hamiltonian for a quasi-Hamiltonian
representation of the system with respect to a non-standard Poisson
structure, coming from the special conformal Killing tensor (scKt
for short). It must be said that scKts for general metrics had
already appeared also, with slightly different assumptions and under
different names, in the work of Benenti (see the comprehensive
review  paper \cite{Benenti2}).

As our title reflects, there is an extra aspect about the cofactor
systems we want to study. So called `driven cofactor systems' again
were introduced by Lundmark and Rauch-Wojciechowski \cite{LW}, still
in the context of mechanical systems with Euclidean kinetic energy
metric and hence having no terms quadratic in the velocities in the
second-order equations of motion. Briefly, the systems discussed in
\cite{LW} are of the form
\begin{eqnarray*}
\ddot{y}^i &=& Q^i(y^j), \qquad i=1,\ldots, m \\
\ddot{x}^a &=& Q^a(y^i,x^b) \qquad a=1,\ldots, n.
\end{eqnarray*}
To begin with, they exhibit a given partial decoupling whereby the
$y$-equations are referred to as the driving system, and the
$x$-equations as the driven part. In addition, it is assumed that
the overall system is of cofactor type and that, more restrictively,
the force terms $Q^a$ come from a potential, parametrically
depending on the driving coordinates. These are rather strong
conditions indeed, but they were shown to lead to quite striking
conclusions in \cite{LW}. First of all the driving system turns out
to be of cofactor type in its own right. Secondly, the driven
system, when regarded as a time-dependent system along solutions
$y(t)$ of the driving system, has $n$ (time-dependent) quadratic
integrals. Most astonishingly, however, the authors managed to show
that (under some technical assumptions) there exists a
time-dependent standard canonical transformation, which has the
effect of shifting the time-dependence to an overall factor, so that
an autonomous Hamiltonian can be identified which turns out to be of
St\"ackel type.

Clearly, such results lead to a double challenge. The first question
is whether they can be extended to general cofactor systems and can
be understood in more intrinsic geometric terms this way. Secondly,
if there exists an intrinsic scheme behind these observations,
coordinates with respect to which the system partially decouples
should not be regarded as part of the data; their existence and a
constructive procedure to find them should follow from testing
coordinate free conditions. This second aspect fortunately has
sufficiently been explored in the literature. Second-order equations
with such a decoupling property were called submersive in \cite{KT},
where both local and global criteria were developed for their
characterization. The local conditions for existence of appropriate
coordinates and for their construction were turned into a more
compact and transparent form by using the geometric calculus along
the tangent bundle projection developed in \cite{MaCaSaIII} (see
also \cite{Martinez}).

In a brief communication \cite{SarBier}, one of us has presented the
main ingredients towards the resolution of the above double
challenge. First of all, \cite{SarBier} provides a coordinate free
definition of a driven cofactor system. It further shows that the
driving system carries a cofactor structure in its own right and
that the overall system has a second, in some sense degenerate
cofactor representation, from which one can develop a scheme to
construct $n+1$ quadratic integrals. One of these is the integral
corresponding to the cofactor nature of the driving system; the
other $n$ (the dimension of the driven system) are integrals of the
driven system along solutions of the driving one. What was not
sufficiently understood at that time was the delicate issue of the
nature of a canonical transformation (if any) which would have the
effect of eliminating the time-dependence form the driven system in
a way which leads to the identification of a St\"ackel system, as in
the Euclidean case. The main purpose of the present paper is
precisely to complete this part of the story. To understand
precisely what is happening in the process, we need much more
explicit information about the recursive procedure which leads to
the $n+1$ first integrals referred to before, and also about the
precise structure of these integrals. The details of such
computations will take a substantial part of our discussion, but in
doing so, the present paper will also complement a number of the
results already mentioned in \cite{SarBier}, where the focus was
more on existence issues.

The next section contains the basic definitions needed to describe
driven cofactor systems and a summary of results presented in
\cite{SarBier}. In Section~3, we develop the algorithm which leads
to the identification of the $n+1$ first integrals referred to
before. A key issue for understanding the nature of the driven
system is the identification of a scKt for its proper metric. It is
shown in Section~4 that this scKt does not, however, turn the driven
system into a cofactor system in its own right; instead it does that
for a modified driven system which turns out to play a role further
on. In Section~5, we start by identifying Darboux coordinates for
the symplectic form associated to the scKt of the complete system,
they are obtained by suitably modifying the momenta. But we
gradually develop arguments then to come to an even better selection
of modified momenta, which takes the specific decoupling properties
of our system into account and are shown to be related to a
time-dependent (standard) canonical transformation for the driven
part of the system. In Section~6, we prove that the application of
this canonical transformation, followed by one which comes from
using eigenfunctions as new coordinates, produces the rather
miraculous effect of reducing the driven system essentially to an
autonomous St\"ackel type system. The proofs in Section~6 are partly
based on simple, indirect arguments, but they are supported also by
explicit computations about the structure of all first integrals,
which are presented in an Appendix. A couple of illustrative
examples are presented in Section~7.

\section{Preliminaries}

The mechanical systems we are talking about in this paper belong to
the class of nonconservative Lagrangian systems, governed by
equations of the form
\begin{equation}
\frac{d}{dt}\left(\fpd{T}{\dot{q}^\alpha}\right) - \fpd{T}{q^\alpha}
= Q_\alpha. \label{nonconservative}
\end{equation}
Here $T=\onehalf g_{\alpha\beta}(q)v^\alpha v^\beta$ is a kinetic
energy function on the tangent bundle $TM$ of a Riemannian manifold
$M$ with metric $g$; the nonconservative forces $Q_\alpha$ are
assumed to depend on the position variables only and thus can be
viewed as components of a 1-form $\mu=Q_\alpha(q)dq^\alpha$ on $M$.
We recall from \cite{CST} that a special conformal Killing tensor
$J$ (scKt) is a type $(1,1)$ tensor field on $M$, which is symmetric
with respect to $g$ and satisfies (lowering an index in the usual
way)
\begin{equation}
J_{\alpha\beta|\gamma} = \onehalf (\alpha_\alpha g_{\beta\gamma} +
\alpha_\beta g_{\alpha\gamma}), \quad \mbox{which further implies
that\ } \alpha = d(\tr J). \label{scKt}
\end{equation}
The Nijenhuis torsion ${\mathcal N}_J$ vanishes, implying that
${d_J}^2=0$. Recall for completeness (see \cite{FroNij}) that $d_J$
is the derivation of degree 1 which (anti)commutes with the exterior
derivative and whose action on functions $f\in\cinfty{M}$ is given
by $d_Jf= J(df)$. The Nijenhuis torsion ${\mathcal N}_J$ is the
vector-valued 2-form with components
\[
J^\alpha_\beta\left(\fpd{J^\beta_\gamma}{q^\delta}-\fpd{J^\beta_\delta}{q^\gamma}\right)-
J^\beta_\delta\fpd{J^\alpha_\gamma}{q^\beta} +
J^\beta_\gamma\fpd{J^\alpha_\delta}{q^\beta}.
\]
We also have that ${D_J}^2=0$, where $D_J$ is the `gauged
differential operator', defined by
\begin{equation}
D_J \rho= d_J \rho+ d(\tr J)\wedge\rho = (\det J)^{-1}d_J((\det
J)\rho). \label{DJ}
\end{equation}
The equality of both expressions in the defining relation of $D_J$
follows from the fact that
\begin{equation}
d_J(\det J) = (\det J) d(\tr J), \label{dJ(detJ)}
\end{equation}
for any tensor with vanishing Nijenhuis torsion.

The concept of a cofactor system on a general Riemannian manifold
was introduced in \cite{CS2001} in the following way.

\begin{dfn} A cofactor system is a triple $(g,\mu,J)$ on a Riemannian manifold where
$g$ is the metric, $\mu$ is a 1-form and $J$ is a nonsingular
special conformal Killing tensor such that $D_J\mu=0$.
\end{dfn}

Combining this concept with the geometric notion of submersiveness,
we came in \cite{SarBier} to the following generalization and
coordinate free formulation of the kind of driven cofactor systems
introduced in \cite{LW}.

\begin{dfn}\ A driven cofactor system is a cofactor system $(g,\mu,J)$,
for which there exists a distribution $K$ along the projection
$\tau: TM\rightarrow M$, with the properties
\begin{equation}
\Phi(K)\subset K, \quad \nabla K \subset K, \quad \DV{Z}K\subset K
\quad\forall Z\in\tvectorfields, \label{K}
\end{equation}
\begin{equation}
d\mu(K,K)=0, \quad \DH{}\!\mu(K^\perp,K)\neq 0. \label{mu}
\end{equation}
\end{dfn}

To understand the meaning of these conditions, one needs to know
about the following intrinsic geometrical concepts associated to
general second-order differential equation fields (\sode s), say
\[
\Gamma=v^\alpha \,\fpd{}{q^\alpha} + f^\alpha(q,v)\,
\fpd{}{v^\alpha}.
\]
Such a \sode\  comes with a connection, determined by the following
horizontal lift construction from $\vectorfields{M}$ (the module of
vector fields on $M$) to $\vectorfields{TM}$
\[
X= X^\alpha(q)\fpd{}{q^\alpha} \mapsto \H{X}= X^\alpha\,H_\alpha,
\]
where
\[
H_\alpha=\fpd{}{q^\alpha} - \Gamma^\beta_\alpha\fpd{}{v^\beta},
\quad \mbox{with}\quad
 \Gamma^\alpha_\beta = - \frac{1}{2}\;\fpd{f^\alpha}{v^\beta}.
\]
$\Gamma^\alpha_\beta$ are the connection coefficients. Naturally, we
also have a vertical lift
\[
X= X^\alpha(q)\fpd{}{q^\alpha} \mapsto \V{X}= X^\alpha\,V_\alpha,
\quad\mbox{where\ } V_\alpha=\fpd{}{v^\alpha}.
\]
It is clear that both operations still make sense if we allow the
components of $X$ to be functions $X^\alpha(q,v)$, meaning that we
extend the domain of the horizontal and vertical lift to
$\tvectorfields$, the $\cinfty{TM}$-module of vector fields along
$\tau:TM\rightarrow M$. These lifts further give rise to
corresponding {\sl horizontal and vertical covariant derivative
operators\/} $\DH{X}$ and $\DV{X}$, determined by the following
action on functions $F\in\cinfty{TM}$ and basic vector fields (and
then further extended by duality):
\begin{eqnarray*}
\DV{X} F = X^\alpha\,V_\alpha(F), &\qquad& \DV{X}\fpd{}{q^\alpha} = 0 \\[1mm]
\DH{X} F = X^\alpha\,H_\alpha(F), &\qquad& \DH{X}\fpd{}{q^\alpha} =
X^\beta V_\alpha(\Gamma^\gamma_\beta) \fpd{}{q^\gamma}\,.
\end{eqnarray*}
The $\DH{}$ in (\ref{mu}) is a covariant differential, defined for
any tensor field ${\mathcal T}$ along $\tau$ by $\DH{}{\mathcal
T}(X,Y,\ldots) = \DH{X}{\mathcal T}(Y,\ldots)$. Furthermore, the
decomposition of $\lie{\Gamma}\H{X}$ into its horizontal and
vertical part identifies the important concepts of {\it dynamical
covariant derivative\/} $\nabla$, a self-dual degree 0 derivation on
tensor fields along $\tau$, and {\it Jacobi endomorphism\/}, a
$(1,1)$ tensor $\Phi$ along $\tau$:
\[
\lie{\Gamma}\H{X} = \H{(\nabla X)} + \V{\Phi(X)} .
\]
For practical purposes, it suffices to know that:
\[
\del F = \Gamma(F) \qquad \del \fpd{}{q^\alpha}=
\Gamma^\beta_\alpha\fpd{}{q^\beta} \qquad \del dq^\alpha = -
\Gamma^\alpha_\beta dq^\beta\,,
\]
\[
\Phi^\alpha_\beta = -\fpd{f^\alpha}{q^\beta} -
\Gamma^\alpha_\gamma\Gamma^\gamma_\beta -
\Gamma(\Gamma^\alpha_\beta)\,.
\]
For the broader picture of derivations of forms along $\tau$, where
the above mentioned concepts play a distinctive role, one can
consult \cite{MaCaSaI, MaCaSaII}.

As was discussed in Proposition~1 of \cite{SarBier} and before that
already in \cite{Martinez}, the existence of a distribution $K$
along $\tau$ having the properties (\ref{K}) implies that $K$ is
actually spanned by a distribution on $M$ which is integrable, and
introducing adapted coordinates $(y^i,x^a)$ has the effect that the
given \sode\ partially decouples into equations of the form
\begin{eqnarray}
\ddot{y}^i &=& f^i(y,\dot{y}), \quad\qquad i=1,\ldots, m, \label{driving} \\
\ddot{x}^a &=& f^a(x,y,\dot{x},\dot{y}), \quad a=1,\ldots, n \ \
(\mbox{here\ }n+m=N=\dim M). \label{driven}
\end{eqnarray}
We will keep referring to the decoupled $y$-system as the {\sl
driving equations\/} and to the remaining $x$-equations as the {\sl
driven system\/}. But there is more to it in this context, which
brings us to the final requirements (\ref{mu}) in definition~2.
$K^\perp$ of course is the orthogonal complement of $K$ with respect
to the Riemannian metric $g$. An important property of $g$ is that
\begin{equation}
\del g =0, \label{delg}
\end{equation}
which follows from the fact that $g$ is the Hessian of the
Lagrangian for the conservative part of the system and the knowledge
that the nonconservative forces determined by $\mu$ have no effect
on the connection coefficients. It was shown in \cite{SarBier} that
as a result of (\ref{delg}) $K^\perp$ inherits the properties
$\nabla K^\perp \subset K^\perp$ and $\DV{Z}K^\perp\subset K^\perp$
from $K$ and that this is enough to conclude that the two
complementary distributions are in fact simultaneously integrable.
Hence we can choose $x^a$ and $y^i$ coordinates which are adapted to
$K$ and $K^\perp$ at the same time, i.e.\ in such a way that
\[
K = {\rm sp\,}\left\{\fpd{}{x^a}\right\}, \qquad K^\perp = {\rm
sp\,}\left\{\fpd{}{y^i}\right\}\,.
\]
It also follows that the kinetic energy part in the equations of
motion (\ref{nonconservative}) decouples completely. In other words,
if we put $g_1:= {\left. g \right|}_{K^\perp}$ and $g_2:= {\left. g
\right|}_{K}$, then in adapted coordinates:
\begin{equation}
g_1 = g_{ij}(y) dy^i\otimes dy^j, \qquad g_2 = g_{ab}(x) dx^a\otimes
dx^b,  \label{gparts}
\end{equation}
while $g_{ia}=g_{ai}=0$. Similarly, for the corresponding connection
coefficients, we have
\[
\Gamma^i_{jk}=\Gamma^i_{jk}(y),\qquad \Gamma^a_{bc}=\Gamma^a_{bc}(x)
\]
and all other connection coefficients are zero. If we insist on
keeping some partial coupling between a driving and a driven part of
the dynamics therefore, this can only come from the nonconservative
forces in $\mu$. The second of the conditions (\ref{mu}) exactly
guarantees such a coupling. The condition $d\mu(K,K)=0$ on the other
hand can easily be seen to model the additional assumption that, in
adapted coordinates, the driven part has force terms $Q_a$ which are
derivable from a potential energy function (with parametric
dependence on the driving coordinates $y^i$).

\section{The cofactor pair scheme on $M$ and $n+1$ quadratic integrals}

It is appropriate to introduce complementary projection operators
\[
P_1: \vectorfields{M} \rightarrow K^\perp, \qquad P_2:
\vectorfields{M} \rightarrow K.
\]
We thus have $P_1+P_2=I_N$ (the identity tensor on the
$N$-dimensional manifold $M$) , $P_1\circ P_2= P_2\circ P_1=0$,
$P_i^2=P_i$, and we occasionally put $P_1|_{K^\perp} = I_m, \
P_2|_K=I_n$. As in \cite{SarBier}, we then look at the scKt tensor
$J$ (for its action on vector fields) as the sum of the following
four parts
\begin{equation}
J_i= P_i\circ J\circ P_i, \ i=1,2, \quad J_{12}=P_1\circ J\circ P_2,
\quad J_{21}=P_2\circ J\circ P_1, \label{Jparts}
\end{equation}
and we shall use a similar notation also for other type (1,1) tensor
fields of interest further on.

{\sc Important notational convention}: In principle, such tensor
fields act on the whole module of vector fields on $M$; we shall use
the same notation, however, when we consider their restriction to
the appropriate submodule $K$ or $K^\perp$ where they are not zero.

Recall that the scKt conditions (\ref{scKt}), when expressed in
terms of the coefficients of the original type $(1,1)$ tensor, take
the form
\begin{equation}
J^\alpha_{\beta|\gamma} := \fpd{J^\alpha_\beta}{q^\gamma} -
J^\alpha_\sigma\Gamma^\sigma_{\beta\gamma} +
J^\sigma_\beta\Gamma^\alpha_{\sigma\gamma} = \onehalf
(\alpha_\beta\delta^\alpha_\gamma + \alpha_\sigma
g^{\sigma\alpha}g_{\beta\gamma}). \label{scKt11}
\end{equation}
Taking the decoupling properties of $g$ into account, it follows
that in adapted coordinates: $\partial J^i_j/\partial x^a= \partial
J^a_b/\partial y^i=0$ and
\begin{equation}
J^a_{i|k} = \fpd{J^a_i}{y^k} - J^a_j\Gamma^j_{ik} = \onehalf
\alpha_b g^{ba}g_{ik}, \qquad J^i_{b|a} = \fpd{J^i_b}{x^a} -
J^i_c\Gamma^c_{ba} = \onehalf \alpha_j g^{ji}g_{ab}. \label{Jprops}
\end{equation}
Hence, the different blocks of $J$ have the following type of
restricted dependence on the adapted coordinates:
\[
J_1 = J^i_j(y) \fpd{}{y^i}\otimes dy^j, \quad J_2 = J^a_b(x)
\fpd{}{x^a}\otimes dx^b,
\]
while
\[
J_{21} = J^a_i(y,x) \fpd{}{x^a}\otimes dy^i, \quad J_{12} =
J^i_a(y,x) \fpd{}{y^i}\otimes dx^a,
\]
with
\begin{equation}
\fpd{J^a_i}{y^k} = \fpd{J^a_k}{y^i} \quad\mbox{and}\quad
\fpd{J^i_a}{x^b} = \fpd{J^i_b}{x^a}. \label{Jprop}
\end{equation}
Furthermore, with $\mu_i:=P_i(\mu)$ and assuming that $J_1$ is
nonsingular, we know from \cite{SarBier} that $(g_1,\mu_1,J_1)$
provides a cofactor representation of the driving system. This
implies in particular (see further for more details) that this
system has a quadratic integral
\begin{equation}
E^1 = \onehalf A^1_{ij}(y)\dot{y}^i\dot{y}^j + W^1(y), \label{E1}
\end{equation}
where $A^1= \cof J_1$. For clarity: the cofactor tensor $A$ of a
type $(1,1)$ tensor $J$ (notation $A=\cof J$) is defined by the
relation $JA=AJ=(\det J)I$. The function $W^1$ is determined by the
relation $A^1\mu_1= - dW^1$, which is locally equivalent to the
condition $D_{J_1}\mu_1=0$ in the definition of a cofactor system.
Incidentally, the minus sign in the expression for $A^1\mu_1$ is in
correspondence with the $+W^1$ in the integral $E^1$ (and
\cite{SarBier} contains a sign error in this sense).

We finally recall from \cite{SarBier} that $P_2$, being a degenerate
type $(1,1)$ tensor field on $M$, formally satisfies the
requirements for a scKt with respect to $g$, and additionally has
the property $D_{P_2}\mu=0$, meaning that we have a second
(degenerate) cofactor representation for the full nonconservative
system $\Gamma$, the implications of which we will investigate now.
The starting point is that $J+aP_2$ also satisfies the scKt
condition for any constant $a$ (and is nonsingular for sufficiently
small values of $a$). Let $A(a)$ denote the cofactor tensor of
$J+aP_2$, so that
\begin{equation}
(J + aP_2) A(a) = A(a)(J + aP_2) =\det(J+aP_2) I_N. \label{cofA(a)}
\end{equation}
Since $P_2=I_n$ in adapted coordinates, it is clear that $A(a)$ and
$\det(J+aP_2)$ are polynomials in $a$ of degree $n$. We represent
them as follows,
\[
A(a) = \sum_{i=1}^{n+1} A_{(i)}a^{i-1}, \qquad
\det(J+aP_2)=\sum_{i=1}^{n+1} \Delta_{(i)}a^{i-1},
\]
and identify the coefficients of equal powers of $a$ in the
identities (\ref{cofA(a)}). We get
\begin{eqnarray}
P_2A_{(n+1)} &\!=\!& A_{(n+1)}P_2 \,=\, 0 , \label{cfn+1} \\
JA_{(i+1)}+ P_2A_{(i)} &\!=\!& A_{(i+1)}J+ A_{(i)}P_2 \,=\,
\Delta_{(i+1)}I_N,
\qquad (1 \leq i \leq n), \label{cfn} \\
JA_{(1)} &\!=\!& A_{(1)}J \,=\, \Delta_{(1)}I_N. \label{cf1}
\end{eqnarray}
Information about the block structure of the different $A_{(i)}$
should follow by left and right actions of the projectors $P_k$ on
these relations. It immediately follows from (\ref{cfn+1}) for
example that
\begin{equation}
A_{(n+1)21} = A_{(n+1)12} = A_{(n+1)2} = 0, \label{An+1}
\end{equation}
while (\ref{cfn}) with $i=n$ subsequently implies that $J_1
A_{(n+1)1} = A_{(n+1)1}J_1 = \Delta_{(n+1)} P_1$, or in the
restriction to $K^\perp$:
\[
J_1 A_{(n+1)1} = A_{(n+1)1}J_1 = \Delta_{(n+1)} I_m.
\]
Taking into account that terms of degree $n$ in $A(a)$ can only be
produced by minors of the $J_1$ elements of $J$, this implies that
\begin{equation}
A_{(n+1)1}=A^1=\cof J_1 \qquad \mbox{and} \qquad \Delta_{(n+1)}=\det
J_1. \label{An+11}
\end{equation}
The same equation (\ref{cfn}) for $i=n$ then further yields
information about parts of $A_{(n)}$: it is indeed not hard to show
by appropriate actions of the projectors that (in the restriction to
$K$ or $K^\perp$)
\begin{equation}
A_{(n)2}=(\det J_1) I_n, \qquad A_{(n)21}= - J_{21} (\cof J_1),
\qquad A_{(n)12}= - (\cof J_1)J_{12}. \label{An}
\end{equation}
For the remaining part of $A_{(n)}$, we have to move to the next
line in the hierarchy ($i=n-1$ in (\ref{cfn})): a right and left
action of $P_1$ leads to $J_1 A_{(n)1} + J_{12}A_{(n)21} =
\Delta_{(n)} P_1$, from which it readily follows that
\begin{equation}
A_{(n)1} = (\det J_1) J_1^{-1}J_{12}J_{21}J_1^{-1} + \Delta_{(n)}
J_1^{-1}. \label{An1}
\end{equation}
We will come back in more detail to the continuation of this
recursive scheme when we are in a position to gather information
about the functions $\Delta_{(i)}$. But note that the final relation
(\ref{cf1}) expresses that $A_{(1)}=\cof J$ and $\Delta_{(1)}=\det
J$.

The important feature about having a cofactor pair scheme is that it
gives rise to a family of quadratic integrals which are in
involution with respect to a double Poisson structure. We sketch how
this will work here by following the procedure explained in
\cite{CS2001}. Let $\hat{g}:TM\rightarrow T^*M$ denote the Legendre
type diffeomorphism coming from the given metric $g$, which in
(general) coordinates reads $p_\alpha= g_{\alpha \beta}v^\beta$, and
put $\hat{\Gamma}= \hat{g}_*\Gamma$, where $\Gamma$ is the \sode\
defined by the cofactor system $(g,\mu,J)$. Then, according to
theorem~2 in \cite{CS2001}, the requirements about a cofactor system
are equivalent to $\hat{\Gamma}$ representing a quasi-Hamiltonian
system with a quadratic Hamiltonian with respect to the Poisson
tensor which can be constructed from the tensor $J$ with vanishing
Nijenhuis torsion. More precisely, if $\tilde{J}$ denotes the
so-called complete lift of $J$ to $T^*M$ and $P_0$ is the canonical
Poisson map on $T^*M$, i.e. the canonical Poisson tensor,
interpreted as a map $P_0: \oneforms{T^*M} \rightarrow
\vectorfields{T^*M}$, we put $P_J=\tilde{J}\circ P_0$ and the point
about $\hat{\Gamma}$ is that
\begin{equation}
(\det J)\hat{\Gamma} = P_J(dH), \quad \mbox{with}\quad H= \onehalf
A^{\alpha\beta}p_\alpha p_\beta +W, \label{quasiHam}
\end{equation}
where $A$ is the cofactor tensor of $J$ and $W$ is the function,
determined (up to a constant) by the property $A\mu=-dW$. Having a
double cofactor representation with special conformal Killing
tensors such as $J$ and $P_2$ in our present situation, will entail
that $\hat{\Gamma}$ satisfies a relation of the form
\begin{equation}
(\det (J+aP_2))\hat{\Gamma} = P_{J+aP_2}dH(a), \quad\mbox{with}\quad
H(a)=\onehalf  A^{\alpha\beta}(a)p_\alpha p_\beta +W(a).
\label{quasiHampara}
\end{equation}
Since $P_{J+aP_2}= P_J + aP_{P_2}$, one can see that $H(a)$ will be
a polynomial of degree at most $n$ in the parameter $a$, say of the
form
\begin{equation}
H(a)= \sum_{i=1}^{n+1} H_{(i)}a^{i-1}, \quad\mbox{with}\quad
H_{(i)}=\onehalf  A^{\alpha\beta}_{(i)}p_\alpha p_\beta +W_{(i)},
\label{H(a)}
\end{equation}
where $A^{\alpha\beta}_{(i)}$ comes from the $A_{(i)}$ tensor
considered before, with an index raised by the given metric $g$, and
$W_{(i)}$ is determined by the relation $A_{(i)}\mu=- dW_{(i)}$.
Naturally $H(a)$ and therefore all its coefficients $H_{(i)}$ will
be first integrals of the system $\hat{\Gamma}$.

Before proceeding, observe that $H_{(n+1)}$ is essentially the
function $E^1$ which is a first integral of the driving system (more
precisely $H_{(n+1)}=g_*E^1$). Indeed, it follows from the block
structure of $A_{(n+1)}$ determined before that in adapted
coordinates, with $p_i=g_{ij}(y)\,\dot{y}^j$,
\begin{equation}
H_{(n+1)}= \onehalf {(A^1)}^{ij}(y)p_i p_j +W^1(y). \label{H(n+1)}
\end{equation}
It is worth mentioning also that for a general $J$, the coordinate
expression of $P_J$ is given by
\begin{equation}
P_J = J^\alpha_\beta \fpd{}{p_\beta}\wedge \fpd{}{q^\alpha} -
\onehalf p_\gamma\left(\fpd{J^\gamma_\alpha}{q^\beta}-
\fpd{J^\gamma_\beta}{q^\alpha}\right)
\fpd{}{p_\alpha}\wedge\fpd{}{p_\beta}. \label{PJ}
\end{equation}
Naturally, in adapted coordinates,
\begin{equation}
P_{P_2}= \fpd{}{p_a}\wedge\fpd{}{x^a} , \label{P2}
\end{equation}
where the momentum variables $p_a$ are defined by $p_a=
g_{ab}(x)\dot{x}^b$.

Taking the various polynomial representations into account, property
(\ref{quasiHampara}) becomes
\[
({\textstyle\sum_{i=1}^{n+1}} \Delta_{(i)}a^{i-1}) \hat{\Gamma} =
(P_J + a P_{P_2})({\textstyle\sum_{i=1}^{n+1}} dH_{(i)}a^{i-1}).
\]
Identifying coefficients of equal powers of $a$ requires first of
all that we should have $P_{P_2}(dH_{(n+1)})=0$, and this is clearly
verified in view of the preceding observations. We further must have
that
\begin{equation}
\Delta_{(i)}\hat{\Gamma} = P_J(dH_{(i)}) + P_{P_2}(dH_{(i-1)}),
\qquad 1<i\leq n+1, \label{Hrecursion}
\end{equation}
and finally for $i=1$ that $\Delta_{(1)}\hat{\Gamma} =
P_J(dH_{(1)})$, but this is merely a confirmation of the
quasi-Hamiltonian structure coming from $J$, since $\Delta_{(1)}=
\det J$ and $H_{(1)}=H$. There is no reason to expect that the $n$
first integrals $H_{(1)}$ up to $H_{(n)}$ of $\hat{\Gamma}$ would
depend on the coordinates $(x^a,p_a)$ only. Nevertheless, they will
be first integrals of the driven system along solutions
$(y^i(t),p_i(t))$ of the driving system. As for the question of
involutiveness, if we adopt the notational convention that $P_J(df)
= \{f,\,\cdot\,\}_J$, it follows from (\ref{Hrecursion}) that
\begin{equation}
\{H_{(i)},H_{(l)}\}_J + \{H_{(i-1)},H_{(l)}\}_{P_2} =
\Delta_{(i)}\hat{\Gamma}(H_{(l)})=0, \qquad 1<i\leq n+1, 1\leq l
\leq n+1, \label{involrecursion}
\end{equation}
and from the two other observations about $H_{(n+1)}$ and $H_{(1)}$
that
\begin{equation}
\{H_{(1)},H_{(l)}\}_J = \{H_{(n+1)},H_{(l)}\}_{P_2} = 0 \qquad
\forall l. \label{startinvolrecursion}
\end{equation}
Using (\ref{startinvolrecursion}), it further follows from
(\ref{involrecursion}) with $l=1$ and $l=n+1$ respectively, that
also
\[
\{H_{(l)},H_{(1)}\}_{P_2} = \{H_{(l)},H_{(n+1)}\}_J =0,
\]
and then a simple recursive argument in (\ref{involrecursion})
finally implies that all $H_{(i)}$ are in involution with respect to
both the $J$-bracket and the $P_2$-bracket. It is worth observing
(see (\ref{P2})) that in adapted coordinates the $P_2$-bracket
formally looks like the standard Poisson bracket in the $(x^a,p_a)$
coordinates so that we have, along solutions of the driving system,
$n$ integrals for the driven system which are in involution in the
standard sense.

\section{An induced special conformal Killing tensor for the
metric of the driven system}

Part of our basic assumptions so far is that both $J$ and $J_1$ are
nonsingular. We shall see now that this naturally leads to the
introduction of another nonsingular type $(1,1)$ tensor $\bar{J}_2$,
which is a kind of deformation of $J_2$. When talking about
nonsingularity here, don't forget the notational convention
specified before! To say that $J_1$ is nonsingular of course only
makes sense when we mean $J_1|_{K^\perp}$; likewise, nonsingularity
of $\bar{J}_2$ will refer to $\bar{J}_2|_K$.

For later reference we look at the adjoint action of $J$ on 1-forms
for a moment. In doing so (as we already tacitly did for the
cofactor tensor $A$) we use the same notation again, i.e.\ do not
write $J^*$ as is sometimes customary. One has to keep in mind,
however, that when compositions are involved (as in the definition
of $J_{12}$ and $J_{21}$), the order of such compositions has to be
reversed. Let $\alpha$ be an arbitrary 1-form on $M$ and put
$\beta=J\alpha=J(P_1\alpha + P_2\alpha)$. To solve such a relation
for $\alpha$ in terms of $\beta$, it is natural to actually compute
$P_1\alpha$ and $P_2\alpha$. We have
\begin{eqnarray*}
P_1\beta &=& J_1(P_1\alpha) + J_{21}(P_2\alpha), \\
P_2\beta &=& J_{12}(P_1\alpha) + J_{2}(P_2\alpha).
\end{eqnarray*}
Since $J_1$ is nonsingular, it follows from the first relation that
\begin{equation}
P_1\alpha = J_1^{-1}(P_1\beta) - J_1^{-1}J_{21}(P_2\alpha),
\label{P1a}
\end{equation}
and substitution of this result in the second relation leads to
\begin{equation}
{\bar J}_2(P_2\alpha) = (P_2 - J_{12}J_1^{-1}P_1)\beta, \label{P2a}
\end{equation}
where ${\bar J}_2$ is defined (for its action on 1-forms) as
\begin{equation}
{\bar J}_2 = J_2 - J_{12}J_1^{-1}J_{21}. \label{barJ2}
\end{equation}
Obviously $\jb$ vanishes on $K^\perp$, but will be nonsingular on
$K$ so that $P_2\alpha$ now can be obtained from (\ref{P2a}) and
substitution in (\ref{P1a}) subsequently gives us $P_1\alpha$. In
fact, what we are looking at here is the following typical
factorization of a matrix with a block structure, this time written
as representing the action of a type $(1,1)$ tensor on vector
fields:
\begin{equation}
\left(\begin{array}{cc}
J_1 & J_{12}\\
J_{21} & J_2
\end{array}\right)
= \left(\begin{array}{cc}
J_1 & 0\\
J_{21} & 1
\end{array}\right)
\left(\begin{array}{cc}
1 & J_1^{-1}J_{12}\\
0 & J_2-J_{21}J_1^{-1}J_{12}
\end{array}\right).
\label{decompositionJ}
\end{equation}
It follows that $\det J = (\det J_1)(\det \jb)$ so that
nonsingularity of $J$ and $J_1$ implies the same for $\jb$. In
adapted coordinates, $\jb$ has components
\begin{equation}
\jb{}^a_b = J^a_b - J^a_i(J_1^{-1})^i_j J^j_b.
\label{componentsbarJ2}
\end{equation}

\begin{prop}\label{prop-j2b} $\jb$ is a (parameter dependent) special conformal
Killing tensor for $g_2$, and its cofactor tensor is $(\det
J_1)^{-1} A_{2}$ with $A = \cof J$.
\end{prop}
\begin{proof} The proof is a straightforward computation, for which it
will be suitable to work in the adapted $(y^i,x^a)$ coordinates. Let
us lower an index to apply the scKt condition in the covariant form
(\ref{scKt}). Keeping in mind that $g_{ai}=g_{ia}=0$, we have
\[
\jb{}_{cb} := g_{ca}\jb{}^a_b = J_{cb} -
J_{ci}(J_1^{-1})^{ij}J_{jb}.
\]
Since $J$ is symmetric, the same is true for $\jb$. Note that $\jb$
depends on both sets of coordinates, but the $y^i$ are regarded as
external parameters for our present considerations. Since $J_1$
depends on the $y^i$ only and the same of course holds for its
inverse (or its cofactor tensor $A_{(n+1)1}$), we get in the first
place that
\[
\jb{}_{cb|a} = J_{cb|a} - J_{ci|a}(J_1^{-1})^{ij}J_{jb} -
J_{ci}(J_1^{-1})^{ij} J_{jb|a}.
\]
It follows from (\ref{scKt}) that
\begin{equation}
J_{cb|a}=\onehalf(\alpha_cg_{ba} + \alpha_bg_{ca}), \quad \mbox{and}
\quad J_{ci|a}= \onehalf \alpha_i g_{ca}. \label{Jci|a}
\end{equation}
We then easily obtain that
\begin{equation}
\jb{}_{cb|a} = \onehalf (\bar{\alpha}_cg_{ba} + \bar{\alpha}_b
g_{ca}), \quad \mbox{with} \quad \bar{\alpha}_c = \alpha_c -
J_{ci}(J_1^{-1})^{ij}\alpha_j, \label{scKt-jb}
\end{equation}
which shows that $\jb$ is a scKt for $g_2$.

If we denote $\cof J$ as in (\ref{quasiHam}) by $A$ (and we have
seen that it is also $A_{(1)}$), we know that $(\det J)g^{ac} =
J^a_\beta A^{\beta c} = J_i^a A^{ic} + J_b^a A^{bc}$ and $0=(\det
J)g^{jc} = J_k^j A^{kc} + J_b^j A^{bc}$. It then follows that
\begin{eqnarray*}
\jb{}^a_b A^{bc}(\det J_1)^{-1} &=& (\det J_1)^{-1} \left(J^a_b
A^{bc}
- J^a_i(J_1^{-1})^i_j J^j_b A^{bc}\right)\\
&=& (\det J_1)^{-1} \left((\det J)g^{ac} - J^a_iA^{ic}
+ J^a_i(J_1^{-1})^i_j J^j_k A^{kc}\right)\\
&=& (\det J_1)^{-1}(\det J) g^{ac}\\
&=& (\det \jb)g^{ac}
\end{eqnarray*}
which proves the last statement of the proposition. \qed
\end{proof}

Note that $\jb$, perhaps rather unexpectedly, does not give rise to
a cofactor system representation of the driven system in a strict
sense. In other words, it does not seem to be true that the forces
$\mu_2$ of the driven system have the property $D_{\jb}\mu_2=0$, or
equivalently that $(\cof \jb) \mu_2$ is closed. At this moment, the
closest we can get to such a property is that an associated driven
cofactor system $(g_2,\bar{\mu}_2,\jb)$ exists, with modified
nonconservative forces $\bar{\mu}_2$.

\begin{prop} A cofactor system, parametrically depending on the coordinates
of the driving system, is determined by the triple
$(g_2,\bar{\mu}_2,\jb)$, where $g_2$ and $\jb$ are as before, and
\begin{equation}
\bar{\mu}_2 = P_2(d W_{(n)}), \label{barmu2}
\end{equation}
$W_{(n)}$ being the function encountered in the recursive scheme
following from (\ref{quasiHampara}).
\end{prop}
\begin{proof} We know that the given nonconservative forces $\mu$
satisfy the relation $A\mu=-dW$, or equivalently $(\det J)\mu =
-J(dW)$. Consider the scheme of the beginning of this section which
led to the introduction of $\jb$. Letting $-dW$ play the role of
$\alpha$ and $(\det J)\mu$ the role of $\beta$, we put
$P_i(dW)=d_iW$ for convenience. The relation (\ref{P2a}) becomes
\begin{equation}
- \jb(d_2W) = (\det J)\mu_2 + (\det \jb)J_{12}(d_1 W^1),
\label{jbd2W}
\end{equation}
where for the second term on the right we have taken into account
that $\mu_1$ satisfies $(\det J_1)\mu_1 = - J_1(d_1W^1)$ in view of
the cofactor representation of the driving system, and that $\det J
= (\det J_1)(\det \jb)$. Secondly, projecting the relation
$A_{(n)}\mu = - dW_{(n)}$ under $P_2$, we get in the first place
that
\[
A_{(n)2}\mu_2 + A_{(n)12}\mu_1 =- d_2W_{(n)}.
\]
Using the information gathered about $A_{(n)}$ in (\ref{An}) and the
cofactor system property of the driving system which was just
recalled, this immediately leads to
\begin{equation}
(\det J_1)\mu_2 + J_{12}(d_1W^1) = - d_2W_{(n)}. \label{d2Wn}
\end{equation}
Substituting this result in (\ref{jbd2W}), we obtain the relation
\begin{equation}
(\det \jb)d_2W_{(n)} = \jb(d_2W), \label{modified-mu}
\end{equation}
which implies, with $\bar{\mu}_2=d_2W_{(n)}$, that
\begin{equation}
(\cof \jb) \bar{\mu}_2 = d_2W. \label{cofactor-barmu}
\end{equation}
Together with the knowledge that $\jb$ is a scKt with respect to
$g_2$, this expresses that the triple $(g_2,\bar{\mu}_2,\jb)$
determines a cofactor system. \qed

\end{proof}

It is worth emphasizing again that the cofactor representation of
this modified driven system must be seen also as a statement about
second-order differential equations for the $x^a$, which
parametrically depend on the $y$-coordinates. This is clear, for
example, from the fact that the exterior derivative in the picture
is $d_2$.

As explained in the introduction, one of our main objectives is to
explore the remarkable situation that the driven system, although
being essentially time-dependent along solutions of the driving
system, does give rise in the end to a St\"ackel-type
Hamilton-Jacobi separability anyway. This will require a
supplementary assumption about existence of independent
eigenfunctions of $J$. We want to show at this point that it is
again the tensor $\jb$ which is relevant for this purpose.

Consider the (degenerate kind of) eigenvalue equation
$\det(J-\lambda P_2)=0$, which is a polynomial equation of degree
$n$ for $\lambda$. Our basic assumption now is that this equation
has $n$ functionally independent solutions $u^a$, in the sense that
the 1-forms $P_2(du^a)$ are linearly independent. If we think of the
$u^a$ as expressed in terms of the adapted $(y^i,x^b)$ coordinates,
this amounts to saying that the Jacobian $(\partial u^a/\partial
x^b)$ is nonsingular. From the identity $\det(J-u^a(y,x) P_2)\equiv
0$ for each fixed $u^a(y,x)$, it follows that
\begin{eqnarray*}
0 &\equiv& d\,(\det(J-u^a(y,x) P_2))\\
&=& d\,(\det(J-\lambda P_2))|_{\lambda=u^a(y,x)} +
\left.\fpd{(\det(J-\lambda P_2))}{\lambda}\right|_{\lambda=u^a(y,x)}
du^a.
\end{eqnarray*}
But since $J-\lambda P_2$ has vanishing Nijenhuis torsion, we know
from (\ref{dJ(detJ)}) that
\[
(J-\lambda P_2)\,d\,(\det(J-\lambda P_2)) = \det(J-\lambda P_2)\,d
\tr(J-\lambda P_2)
\]
for all $\lambda$. It then follows by acting with $(J-u^a(y,x) P_2)$
on the preceding identity that
\[
\left.\fpd{(\det(J-\lambda
P_2))}{\lambda}\right|_{\lambda=u^a(y,x)}\,(J-u^a(y,x) P_2)\,du^a=0,
\]
and thus, since all eigenfunctions are assumed to be simple, that
\begin{equation}
(J-u^a(y,x) P_2)\,du^a=0. \label{eigenformsJ}
\end{equation}

\begin{prop} Assume that the equation $\det(J-\lambda P_2)=0$ has $n$
functionally independent eigenfunctions $u^a$. Then $du^a$ is an
eigenform of $J$ (in the sense of equation (\ref{eigenformsJ}))
corresponding to the eigenvalue $u^a$. Moreover, the $u^a$ are also
eigenfunctions of $\jb$, with $P_2(du^a)$ as corresponding
eigenform.
\end{prop}
\begin{proof} It remains to prove the statement about $\jb$. For
that, it suffices to go back once more to the analysis about
$\beta=J\alpha$ at the beginning of this section, with $du^a$ in the
role of $\alpha$ and $u^a P_2(du^a)$ (no sum!) in the role of
$\beta$. The relation (\ref{P2a}) then says that
\[
\jb (P_2(du^a))=u^a P_2(du^a),
\]
which is precisely what we need. \qed
\end{proof}

\section{A symplectic view-point and Darboux coordinates}

Since $J$ is assumed to be nonsingular, the Poisson tensor
associated to $P_J$ actually comes from a symplectic form which we
call $\omega_J$. The sign convention which we adopt here is that for
any function $F$
\[
X= P_J(dF) \quad \Longleftrightarrow \quad i_X\omega_J = - dF.
\]
One easily verifies that, referring to the general coordinate
expression (\ref{PJ}) of $P_J$, $\omega_J$ is given by
\begin{equation}
\omega_J = {J^{-1}}^\alpha_\beta dp_\alpha\wedge dq^\beta - \onehalf
p_\gamma\left(\fpd{J^\gamma_\alpha}{q^\beta} -
\fpd{J^\gamma_\beta}{q^\alpha}\right){J^{-1}}^\beta_\sigma
{J^{-1}}^\alpha_\rho
 dq^\sigma\wedge dq^\rho. \label{omegaJ}
\end{equation}
So far, this correspondence is valid for any nonsingular type
$(1,1)$ tensor $J$ on $M$. The first term in $\omega_J$ strongly
suggests introducing new momentum variables $\check{p}_\alpha$ by
\[
p_\alpha = J^\beta_\alpha(q) \check{p}_\beta.
\]

\begin{lemma} The coordinate change $(q,p) \leftrightarrow
(q,\check{p})$ determines a Darboux chart for $\omega_J$ if and only
if ${\cal N}_J=0$.
\end{lemma}
\begin{proof}
From $\check{p}_\beta={J^{-1}}^\alpha_\beta p_\alpha$, it is easy to
compute that
\[
d\check{p}_\beta\wedge dq^\beta = {J^{-1}}^\alpha_\beta
dp_\alpha\wedge dq^\beta - \onehalf \check{p}_\delta
\left(\fpd{J^\delta_\tau}{q^\sigma}{J^{-1}}^\tau_\rho -
\fpd{J^\delta_\tau}{q^\rho}{J^{-1}}^\tau_\sigma\right) dq^\sigma
\wedge dq^\rho.
\]
Subtracting this from (\ref{omegaJ}), one easily obtains that the
coefficient of the resulting 2-form (after multiplication by two
$J$-factors) will be zero if and only if ${\cal N}_J=0$. \qed
\end{proof}

As explained in Section~3, the fact that the \sode\ $\Gamma$ on $TM$
satisfies the requirements of a cofactor system is equivalent to
saying that its image $\hat{\Gamma}$ under the Legendre map
$\hat{g}$ has a quasi-Hamiltonian representation (\ref{quasiHam})
with respect to $P_J$. This in turn translates within the symplectic
view-point to $(\det J)i_{\hat{\Gamma}}\omega_J = - dH$. The
preceding lemma then says that $\omega_J$ will take the form of the
standard symplectic form on $T^*M$ when expressed in the variables
$(q,\check{p})$. However, it is not clear that we can take advantage
of such a coordinate change because it does not take account of the
special feature of partial decoupling which our system exhibits. We
shall show that there is a better choice of new momenta, which is
inspired by the above transition to Darboux coordinates but does
take the extra features of a driven cofactor system into account.
Note in passing that there exists a different technique for
obtaining a Hamiltonian representation out of a quasi-Hamiltonian
one: roughly it consists of absorbing the overall factor by a change
of timescale,  which has the disadvantage, however, that the
definition of the new time makes sense only along (the unknown)
solutions of the system. This technique is well documented in
\cite{Benenti2}, for example, and was extensively used in the
context of cofactor systems in \cite{MarBlas}. We believe that the
line of approach we adopt here offers more insight in understanding
the delicate aspects of the driven nature of our cofactor system.

Going back to the \sode\  $\Gamma$ on $TM$, we can first pass to the
coordinates $(y^i,x^a)$ adapted to the complementary distributions
$K^\perp$ and $K$, before passing to the quasi-Hamiltonian
representation $\hat{\Gamma}$. The projectors $P_1$ and $P_2$ have
corresponding actions on $TM$ through their complete lifts; they
give rise to a partial splitting of $\Gamma$ in the form $\Gamma=
\Gamma_1+\Gamma_2$ say, as exhibited in the equations
(\ref{driving}, \ref{driven}). Likewise, the complete lifts
$\tilde{P}_1$, $\tilde{P}_2$ of the projectors to $T^*M$ produce a
partial decoupling $\hat{\Gamma}=\hat{\Gamma}_1+\hat{\Gamma}_2$,
which in adapted coordinates is simply the effect of transforming
(\ref{driving}, \ref{driven}) to equivalent first-order equations by
passing to the momenta $p_i=g_{ij}(y)\dot{y}^j$ and
$p_a=g_{ab}(x)\dot{x}^b$. It is worth illustrating this in more
detail as follows. The \sode\ $\Gamma$ associated to the equations
of motion of the form (\ref{nonconservative}), after applying the
overall Legendre map $\hat{g}$, transforms to
\begin{equation}
\hat{\Gamma}= g^{\alpha\beta}p_\beta \fpd{}{q^\alpha} +
\Gamma^\gamma_{\mu\alpha}g^{\mu\delta}p_\gamma p_\delta
\fpd{}{p_\alpha} + Q_\alpha\fpd{}{p_\alpha}. \label{hatGamma}
\end{equation}
In adapted $(y,x)$-coordinates, in view of the way the components of
$g$ and the connection coefficients decouple (see (\ref{gparts}) and
its consequences), this expression becomes
\begin{eqnarray}
\hat{\Gamma} &=& g^{ij}(y)p_j \fpd{}{y^i} +
\Gamma^k_{ij}(y)g^{il}(y)p_k p_l
\fpd{}{p_j} + Q_j(y)\fpd{}{p_j} \nonumber \\
&& \mbox{} + g^{ab}(x)p_b \fpd{}{x^a} + \Gamma^c_{ab}(x)g^{ad}(x)p_c
p_d \fpd{}{p_b} + Q_b(x,y)\fpd{}{p_b}. \label{hatGyx}
\end{eqnarray}

The first line reflects the fact that the driving system has its own
cofactor representation, i.e.\ satisfies
\begin{equation}
(\det J_1) \hat{\Gamma}_1 = P_{J_1}( dH_{(n+1)}),
\label{quasiHam-driving}
\end{equation}
with $H_{(n+1)}$ as in (\ref{H(n+1)}). Concerning the second line,
we should take into account the extra assumption that the driven
system has a standard Hamiltonian representation: as indicated
before, the condition $d\mu(K,K)=0$ in definition~2 expresses that
the force terms $Q_a$ of the driven system are derivable from a
potential energy function $V(x,y)$ say, depending parametrically on
the driving coordinates $y^i$. It is then clear that the second line
simply expresses that
\begin{equation}
\hat{\Gamma}_2 = P_{P_2}(dh), \quad\mbox{with}\quad h= \onehalf
g^{ab}(x)p_ap_b + V(x,y), \label{Ham-driven}
\end{equation}
keeping in mind that $P_{P_2}$ in adapted coordinates merely is the
standard Poisson structure in the variables $(x^a,p_a)$. With this
splitting of $\hat{\Gamma}$ in mind, it looks more appropriate not
to spoil the decoupled feature of the driving system by introducing
Darboux coordinates for the overall symplectic structure $\omega_J$.
Instead, we can put $p_i= {J_1}^j_i\tilde{p}_j$, which will have the
effect of introducing Darboux coordinates for $\omega_{J_1}$. As for
the driven part $\hat{\Gamma}_2$, let us first investigate in detail
what the introduction of the momenta $\check{p}$ would do.

Formally we can regard the transformation formulas $p_\alpha =
J^\beta_\alpha(q) \check{p}_\beta$ as representing a relation
between 1-forms on $M$, of the type $\beta=J\alpha$ discussed at the
beginning of Section~4. It then follows from the considerations
leading to (\ref{P2a}) that
\begin{equation}
\bar{J}_2{}_a^b \check{p}_b = p_a - J^i_a J_1^{-1}{}^j_i p_j.
\label{tildep}
\end{equation}
This suggests that the more relevant momentum variables for the
driven system actually are $\tilde{p}_a := \bar{J}_2{}_a^b
\check{p}_b$. The conclusion from this preliminary analysis is that
we shall consider the following linear change of momenta
\begin{eqnarray}
p_i &=& {J_1}^j_i\tilde{p}_j, \label{tildep_i} \\
p_a &=& \tilde{p}_a + J^i_a \tilde{p}_i. \label{tildep_a}
\end{eqnarray}
It turns out that the transition from $p_a$ to $\tilde{p}_a$, viewed
as time-dependent transformation along solutions of the driving
system, actually represents a time-dependent canonical
transformation for the driven system in the standard sense and hence
is ideally suited to preserve the special assumption on that system.
Indeed, in view of the second of the properties (\ref{Jprop}), we
know that the components $J^i_a$ of $J$ can be written as $J^i_a=
\partial \psi^i/\partial x^a$ for some functions $\psi^i(x,y)$.
Defining $F(x,\tilde{p},t)$ by
\begin{equation}
F(x,\tilde{p},t) = x^a\tilde{p}_a + \Psi(x,t), \quad\mbox{with}\quad
\Psi(x,t) = \psi^i(x,y(t)) \tilde{p}_i(t), \label{F2}
\end{equation}
we create a generating function of mixed type (depending on the old
position variables $x$ and new momenta $\tilde{p}$) for a standard
canonical transformation $(x^a,p_a)\leftrightarrow
(x^a,\tilde{p}_a)$, which does not change the coordinates,
transforms the momenta according to (\ref{tildep_a}), but must be
viewed as time-dependent along solutions of the driving system. The
Hamiltonian of the transformed system then is given by
\begin{equation}
\tilde{h}(x,\tilde{p},t):= h + \fpd{F}{t}. \label{tildeh}
\end{equation}
We shall see in the next section that this canonical transformation
is one of two steps which are required to relate the original
Hamiltonian $h$ of the driven system to the first integral $H_{(n)}$
of $\hat{\Gamma}$ and that $H_{(n)}$ is the key to understanding the
subtle way in which the driven system in the end corresponds to an
autonomous Hamiltonian system which is separable in the
Hamilton-Jacobi sense.

\section{Separability of the Hamilton-Jacobi equation
for the driven system}

Let us start by computing the function $H_{(n)}$ expressed in the
variables $(y^i,x^a,\tilde{p}_i,\tilde{p}_a)$. We have
\[
H_{(n)} = \onehalf A_{(n)}^{ab} p_ap_b + A_{(n)}^{ai} p_ap_i +
\onehalf A_{(n)}^{ij} p_ip_j + W_{(n)}.
\]
From (\ref{An}), raising an index, we learn that
\[
A_{(n)}^{ab} = (\det J_1) g^{ab}, \qquad A_{(n)}^{ai} = - {A^1}^i_l
J^{la}.
\]
Making the substitutions (\ref{tildep_i}, \ref{tildep_a}) it then
readily follows (remember that $A^1$ is $\cof J_1$) that
\[
H_{(n)} = \onehalf (\det J_1)\,g^{ab}\tilde{p}_a\tilde{p}_b -
\onehalf (\det J_1)J^{bi}J^j_b \tilde{p}_i\tilde{p}_j + \onehalf
A_{(n)}^{kl}{J_1}^i_k {J_1}^j_l\tilde{p}_i\tilde{p}_j  + W_{(n)}.
\]
Using (\ref{An1}) now, we arrive at the following result:
\begin{equation}
H_{(n)} = \onehalf (\det J_1)\,g^{ab}\tilde{p}_a\tilde{p}_b +
\onehalf \Delta_{(n)}J^{ij}\tilde{p}_i\tilde{p}_j + W_{(n)}.
\label{H(n)}
\end{equation}
So, introducing the new variables has the interesting effect of
eliminating the terms in mixed momenta in $H_{(n)}$. Obviously,
however, the effect on $h$ will be the opposite. We get
\begin{equation}
h = \onehalf g^{ab}\tilde{p}_a\tilde{p}_b +
J^{ai}\tilde{p}_a\tilde{p}_i + \onehalf
J^{ai}J^j_a\tilde{p}_i\tilde{p}_j + V. \label{h}
\end{equation}
But for the interpretation as time-dependent canonical
transformation, we need to look at the function $\tilde{h}$ and will
show now that this function is more closely related to $H_{(n)}$.

\begin{lemma}
Under the canonical transformation with generating function
(\ref{F2}), the transformed Hamiltonian $\tilde{h}$ of the driven
system takes the following form, to within an additive function of
time,
\begin{equation}
\tilde{h}= (\det J_1)^{-1} H_{(n)} + J^{ai}\tilde{p}_a\tilde{p}_i.
\label{fulltildeh}
\end{equation}
\end{lemma}
\begin{proof}
We need to add to the expression (\ref{h}) for $h$ the term
\[
\fpd{F}{t}= \fpd{\Psi}{t} = \fpd{\Psi}{y^k}\,\dot{y}^k +
\psi^k\dot{\tilde{p}}_k,
\]
computed along solutions of the driving equations. Since the
$\tilde{p}_i$ were introduced to provide Darboux coordinates, the
Poisson map $P_{J_1}$ in (\ref{quasiHam-driving}) takes the form of
the standard Poisson map so that the equations of the driving system
become
\[
(\det J_1)\, \dot{y}^k = \fpd{H_{(n+1)}}{\tilde{p}_k}, \qquad (\det
J_1)\, \dot{\tilde{p}}_k = - \fpd{H_{(n+1)}}{y^k},
\]
whereby the function $H_{(n+1)}$ from (\ref{H(n+1)}), when expressed
in the new momenta, reads
\begin{equation}
H_{(n+1)} = \onehalf (\det J_1) J_1^{ij} \tilde{p}_i \tilde{p}_j +
W^1. \label{newH(n+1)}
\end{equation}
It is now fairly straightforward to compute that $\tilde{h}$ can be
written as,
\begin{eqnarray}
\tilde{h} &=& h + \fpd{\psi^i}{y^k}J_1^{kj}\tilde{p}_i\tilde{p}_j +
J_1^{j\,l}\Gamma^i_{kl}\psi^k\tilde{p}_i\tilde{p}_j
- \onehalf {J_1^{ij}}_{|k}\psi^k \tilde{p}_i\tilde{p}_j \nonumber \\
&& \hspace{5mm} - \onehalf (\det J_1)^{-1}\fpd{\det J_1}{y^k} \psi^k
J_1^{ij}\tilde{p}_i\tilde{p}_j - (\det J_1)^{-1}\psi^k
\fpd{W^1}{y^k}, \label{tildehex}
\end{eqnarray}
with $h$ as in (\ref{h}). Concerning the function we want to match
it with, we first of all need more info about $\Delta_{(n)}$. Using
the representation (\ref{decompositionJ}) of $J$, we can put
\begin{equation}
J + a P_2 = \left(\begin{array}{cc}
J_1 & 0\\
J_{21} & 1
\end{array}\right)
\left(\begin{array}{cc}
1 & J_1^{-1}J_{12}\\
0 & \jb + a I_n
\end{array}\right),
\label{decompositionJ+aP2}
\end{equation}
from which it follows that
\begin{equation}
\sum_{i=1}^{n+1} \Delta_{(i)}a^{i-1}= \det(J+aP_2)= (\det J_1)
\det(\jb + a I_n), \label{Deltai}
\end{equation}
and hence that $\Delta_{(n)}= (\det J_1) (\tr \jb)$. It follows from
(\ref{H(n)}) that
\begin{equation}
(\det J_1)^{-1}H_{(n)} = \onehalf \,g^{ab}\tilde{p}_a\tilde{p}_b +
\onehalf (\tr \jb)J^{ij}\tilde{p}_i\tilde{p}_j + (\det
J_1)^{-1}W_{(n)}. \label{newH(n)}
\end{equation}
It looks by far not obvious that the expressions (\ref{tildehex})
and (\ref{newH(n)}) would differ only by one term (up to irrelevant
functions of time only). We shall compare them indirectly by
computing their derivatives with respect to the $x^a$ and
$\tilde{p}_a$. The latter is easy and gives
\[
\fpd{\tilde{h}}{\tilde{p}_a} = g^{ab}\tilde{p}_b +
J^{ai}\tilde{p}_i, \qquad \fpd{}{\tilde{p}_a}\left((\det
J_1)^{-1}H_{(n)}\right) = g^{ab}\tilde{p}_b,
\]
which is in line with the result we want to prove. For the other
derivatives, the computations can be written in a somewhat more
compact form if we use the basis of vector fields
\[
X_a = \fpd{}{x^a} + \Gamma^c_{ab}\tilde{p}_c \fpd{}{\tilde{p}_b}
\]
adapted to the connection, rather than the coordinate derivatives
with respect to $x^a$. One can verify, recalling that
$J^k_a=\partial \psi^k/\partial x^a$, that $X_a(\tilde{h})$ can be
written as
\begin{eqnarray*}
X_a(\tilde{h}) &=& J^{bi}_{|a}\tilde{p}_b \tilde{p}_i +
\onehalf\left( J^{bi}_{|a}J^j_b + J^{bi}J^j_{b|a}\right)
\tilde{p}_i \tilde{p}_j  \\
&& \mbox{} + J^i_{a|k}J_1^{jk} \tilde{p}_i\tilde{p}_j - \onehalf
J^k_a {J_1^{ij}}_{|k} \tilde{p}_i\tilde{p}_j - \onehalf (\det
J_1)^{-1} \fpd{\det J_1}{y^k} J^k_a J_1^{ij}
\tilde{p}_i\tilde{p}_j  \\
&& \mbox{} + \fpd{V}{x^a} - (\det J_1)^{-1} J^k_a \fpd{W^1}{y^k}.
\end{eqnarray*}
Making use of the scKt properties of $J$ and $J_1$, plus the
property (\ref{dJ(detJ)}) for $J_1$, this expression considerably
simplifies and finally reduces to
\begin{eqnarray*}
X_a(\tilde{h}) &=& \onehalf \alpha_k g^{ik}\tilde{p}_i\tilde{p}_a -
\onehalf \alpha_l {(J_1^{-1})}^l_k J^k_a J_1^{ij}
\tilde{p}_i\tilde{p}_j
+ \onehalf \alpha_a J_1^{ij} \tilde{p}_i\tilde{p}_j \\
&& \mbox{} + \fpd{V}{x^a} - (\det J_1)^{-1} J^k_a \fpd{W^1}{y^k},
\end{eqnarray*}
where $\alpha$ as before stands for $d(\tr J)$. Note in passing
that, for example, $J^{ij}\equiv J_1^{ij}$.

The computation of $X_a\left((\det J_1)^{-1}H_{(n)}\right)$ is much
easier and gives
\[
X_a\left((\det J_1)^{-1}H_{(n)}\right) = \onehalf \bar{\alpha}_a
J_1^{ij} \tilde{p}_i\tilde{p}_j + (\det J_1)^{-1}
\fpd{W_{(n)}}{x^a},
\]
with $\bar{\alpha} = d(\tr \jb)$. We need to make three more
observations now. The first is that
\[
X_a(J^{bi}\tilde{p}_i\tilde{p}_b) =
J^{bi}_{|a}\tilde{p}_i\tilde{p}_b = \onehalf \alpha_k
g^{ik}\tilde{p}_i\tilde{p}_a,
\]
which takes account of the first term of $X_a(\tilde{h})$. Secondly,
we recall the difference between $\bar{\alpha}_a$ and $\alpha_a$, as
obtained in (\ref{scKt-jb}), which makes that the first term of
$X_a\left((\det J_1)^{-1}H_{(n)}\right)$ matches two terms of
$X_a(\tilde{h})$. Finally, the terms not containing momenta also
match as a result of (\ref{d2Wn}), taking into account that $\mu_2 =
- d_2V$. The conclusion now is that
\[
X_a(\tilde{h}) = X_a\left((\det J_1)^{-1}H_{(n)}+
J^{bi}\tilde{p}_i\tilde{p}_b\right).
\]
Remember that the two functions under consideration are regarded
here as depending on the $x^a$, $\tilde{p}_a$ and time $t$ (along
solutions of the driving equations). Since their derivatives with
respect to the $x^a$ and $\tilde{p}_a$ are the same, the conclusion
is that they are indeed equal up to an additive function of time.
\qed
\end{proof}

The idea now is to try to get rid of the second term on the right in
(\ref{fulltildeh}) by a further, suitable canonical transformation.
It is of some interest to look at this kind of question in all
generality and to observe that it can be resolved indeed by a
suitable point transformation, which necessarily must be
time-dependent however. The notations used in discussing this
general question below have nothing to do with any of the specific
situations encountered so far.

\begin{lemma}
Suppose that $H_1(q,p,t)$ and $H_2(q,p,t)$ are two functions which
differ by terms linear in the $p_i$. Then, there exists a point
transformation, $(q,p) \leftrightarrow (Q,P)$ say, such that the
transformed Hamiltonian of the system with Hamiltonian $H_1$ becomes
the function $H_2$ expressed in the new variables.
\end{lemma}
\begin{proof}
By assumption, we have $H_1 = H_2 + \rho^i(q,t) p_i$ for some
functions $\rho^i$. If $F(q,P,t)=P_iQ^i(q,t)$ is the generating
function of an as yet unspecified point transformation, we know that
the transformed Hamiltonian of the system with Hamiltonian $H_1$
will be given by
\[
\tilde{H}_1 = H_1 + \fpd{F}{t} = H_2 + \rho^i\fpd{Q^j}{q^i}P_j +
\fpd{Q^j}{t}P_j,
\]
so that the desired effect requires that each of the $Q^j$ satisfies
\[
\fpd{Q^j}{t} + \rho^i\fpd{Q^j}{q^i} = 0.
\]
In other words, we need $n$ functionally independent first integrals
of the equations $\dot{q}^i=\rho^i(q,t)$, which can be done in
principle and is of course the same as saying that we have to
integrate those equations. Note that even if the given $\rho^i$
would not depend on time, this procedure can only work with a
time-dependent canonical transformation. \qed
\end{proof}

In the case of interest, we are looking with equation
(\ref{fulltildeh}) at linear terms of the form
\[
\rho^a(x,t)\tilde{p}_a = J^{ai}(x,y(t))\tilde{p}_i(t)\,\tilde{p}_a.
\]
We can try to find first integrals of the equations
$\dot{x}^a=\rho^a(x,t)$ of the form $u^a=u^a(y(t),x)$, i.e.\ which
are such that the time-dependence originates from solutions
$y(t)$ of the driving equations. This means that we have
$\dot{y}^k=J^{ki}\tilde{p}_i$, and the first integral condition
becomes
\begin{equation}
\fpd{u^a}{t}+\rho^b\fpd{u^a}{x^b} =
\tilde{p}_i\left(J^{ik}\fpd{u^a}{y^k} + J^{ib}\fpd{u^a}{x^b}\right)
= 0. \label{firstint}
\end{equation}
We can now prove one of our main results, for which we go back to
the supplementary assumption of the end of Section~4. The
eigenfunctions $u^a(y,x)$ which were introduced there will now be
used as new coordinates for the driven system, along solutions
$y(t)$ of the driving system, and we will denote the corresponding
conjugate momenta (for reference to the notations used in the
Euclidean case in \cite{LW}) by $s_a$.

\begin{prop}
Assume that the equation $\det(J-\lambda P_2)=0$ has $n$
functionally independent solutions $u^a(y,x)$. Then, the canonical
transformation $(x^a,p_a)\leftrightarrow (x^a,\tilde{p}_a)$ with
generating function (\ref{F2}), followed by the canonical
transformation $(x^a,\tilde{p}_a)\leftrightarrow (u^a, s_a)$ with
generating function $F(x,s,t)= s_a u^a(y(t),x)$, has the effect of
transforming the Hamiltonian $h$ of the driven system into the
function $(\det J_1)^{-1} H_{(n)}$.
\end{prop}
\begin{proof}
We know from Lemma~2 that the first step brings the Hamiltonian $h$
into the form (\ref{fulltildeh}). According to Lemma~3, the second
step will eliminate the linear terms in the momenta in
(\ref{fulltildeh}), provided the functions $u^a(y,x)$ have the
property of making the right-hand side of (\ref{firstint}) vanish.
But the $u^a$ satisfy the relations (\ref{eigenformsJ}), from which
it follows by a left action of $P_1$ that
\begin{equation}
J_1 P_1(du^a) + J_{21} P_2(du^a)=0 \quad\mbox{or}\quad
\left(J^k_j\fpd{u^a}{y^k} + J^b_j\fpd{u^a}{x^b}\right)dy^j=0.
\label{auxprop4}
\end{equation}
The desired result now follows by raising an index. \qed
\end{proof}

At this point, it is important to be aware of another general and in
fact very simple result.

\begin{lemma}
Assume that a Hamiltonian system has a Hamiltonian $K(q,p,t)$ of the
form $K=\gamma(t)H(q,p,t)$, whereby $H$ is a first integral of the
system and $\gamma$ is an arbitrary function of $t$ only. Then $H$
in fact cannot explicitly depend on time and the time-dependent
Hamilton-Jacobi equation for $K(q,p,t)$ reduces to the autonomous
one for $H$.
\end{lemma}
\begin{proof}
It is obvious that (with respect to the standard Poisson bracket) we
have $\{K,H\}=\gamma\{H,H\}=0$, so that the property $\dot{H}=0$
reduces to $\partial H/\partial t=0$. The Hamilton-Jacobi equation
for $K$ then reads
\[
\gamma(t)H\left(q,\fpd{S}{q}\right) + \fpd{S}{t}=0,
\]
and looking for a complete solution $S(q,t,\alpha_i)$ of the form
$S=W(q,\alpha_i) - \alpha_1 \int\!\gamma dt$ immediately reduces it
to the Hamilton-Jacobi equation
\[
H\left(q,\fpd{W}{q}\right)= \alpha_1
\]
for the autonomous function $H$. \qed
\end{proof}
Since $\det J_1$, along solutions of the driving system, is a
function of time only and $H_{(n)}$ is known to be a first integral
of the driven system under the same circumstances, it is clear that
the assumptions of Proposition~4 precisely bring us in a situation
where we can draw the quite surprising conclusion that $H_{(n)}$
will no longer be time-dependent in the $(u,s)$ coordinates and
separability of the Hamilton-Jacobi equation for the driven system
is essentially a matter of separability of $H_{(n)}$. It remains to
convince ourselves that the Hamilton-Jacobi equation for a system
with $H_{(n)}$ as Hamiltonian is indeed separable. Now comes a
rather subtle point in the argumentation. The point is this! We want to test
separability of $H_{(n)}$ by using criteria which have an intrinsic,
i.e.\ coordinate independent meaning. As such, it is the function
$H_{(n)}$ which matters, expressed in any kind of canonical
coordinates. Transformation formulas from one set of coordinates to
another can depend on external parameters then, if needed, but
should not be regarded as depending on time because time-dependent canonical
transformations do more to the Hamiltonian function than just
expressing it in the new variables.

To be concrete now, the sufficient conditions for separability which
we want to invoke are intrinsic indeed: they comprise the existence
of a special conformal Killing tensor $J$ for the kinetic energy
metric of the Hamiltonian, plus a corresponding condition for
admissible potentials $V$. The condition $d(JdV)=0$ in
\cite{Benenti97} for example, which actually corresponds to the
particular case of the cofactor condition $A\mu=-dW$ when the forces
are conservative. It is well known (see e.g.\ \cite{Crampin2002})
that if the scKt involved in these conditions has functionally
independent eigenfunctions, then these are separation coordinates
for the Hamilton-Jacobi equation. It so happens that we have needed
these $u^a$-coordinates already to prove in Proposition~4 that the
original Hamiltonian $h$ of the driven system can be transformed
into the function $H_{(n)}$ (up to a factor). But having shown this
way that Hamilton-Jacobi separability becomes a question about the
function $H_{(n)}$, it is more appropriate to put this function to
the test in the set of canonical coordinates $(x^a, \tilde{p}_a)$
which naturally presents itself prior to introducing the separation
coordinates. The delicate issue alluded to above is that, when we
subsequently want to pass to the new variables $(u^a,s_a)$ again,
the interpretation for this part of the story is that the functions
$u^a(y,x)$ are regarded then as depending parametrically on the
$y$-coordinates of the driving system (not as a time-dependent
transformation along solutions $y(t)$ of that system).

Going back to the expression (\ref{newH(n)}) of $H_{(n)}$, taking
into account that the function $H_{(n+1)}$ in (\ref{newH(n+1)}) is
actually a constant parameter (along solutions of the driving
system) which we called $E^1$ in (\ref{E1}), we have that
\begin{equation}
H_{(n)}= \onehalf (\det J_1)\,g^{ab}\tilde{p}_a\tilde{p}_b + (\tr
\jb)(E^1 - W^1(y)) + W_{(n)}(y,x). \label{newerH(n)}
\end{equation}
This is the right expression for activating our test because all the
ingredients we need for that have been prepared in Section~4.

\begin{prop}
$\jb$ is a special conformal Killing tensor for the metric
associated to the quadratic terms in (\ref{newerH(n)}) and the
remaining terms satisfy the conditions for an admissible potential
for Hamilton-Jacobi separability.
\end{prop}
\begin{proof}
We know from Proposition~1 that $\jb$ is a scKt for
$g_2=(g_{ab}(x))$ and recall that the characterizing property
(\ref{scKt}) of scKts, when expressed in terms of the underlying
type $(1,1)$ tensor field is given by (\ref{scKt11}), with
$\alpha=d(\tr J)$. It is then clear that the same $J$ is also a scKt
with respect to any constant multiple of $g$. In the case of the
quadratic terms in (\ref{newerH(n)}), we are precisely looking at a
constant multiple of the metric $g_2$ since $(\det J_1)$ is a
function of the external $y$-parameters only, whence the conclusion
about $\jb$. We further observed above that the condition for
admissible potentials is a reduced form of the cofactor condition
that $(\cof J)\mu$ should be closed. For the situation at hand we
already know from Proposition~2 that the function $W_{(n)}$
satisfies this condition with respect to the scKt $\jb$. It remains
to show that the same is true for the remaining term in
(\ref{newerH(n)}) which in fact, since the factor $(E^1-W^1(y))$ can
be treated as a constant here, amounts to saying that the function
$\tr\jb$ satisfies the condition. But this is trivially the case,
because it follows from (\ref{dJ(detJ)}) that for any tensor $J$
with vanishing Nijenhuis torsion $(\cof J)d(\tr J) = d(\det J)$.
\qed
\end{proof}

We sum up the main results about the driven system now. We know
since Section~3 that the driven system, along solutions of the
driving system, has $n$ integrals $H_{(i)}$ which are in involution
with respect to two Poisson structures, one of which is the standard
one when using coordinates adapted to the integrable distributions
$K$ and $K^\perp$. Under the assumption that the characteristic
equation $\det(J-\lambda P_2)=0$ has $n$ functionally independent
solutions $u^a(y,x)$, we have seen in Proposition~4 that the given
Hamiltonian $h$ of the driven system, can be transformed into the
function $(\det J_1)^{-1} H_{(n)}$. The fact that all the $H_{(i)}$
are first integrals then implies that
\[
\fpd{H_{(i)}}{t} + (\det J_1)^{-1}\{H_{(i)},H_{(n)}\} = 0.
\]
But this, in view of the involutivity, actually means that none of
the integrals will be time-dependent when expressed in the canonical
coordinates $(u^a,s_a)$. This is in line with the results of
Proposition~5 which mainly says that due to the existence of a scKt
$\jb$, the Hamilton-Jacobi equation of the (autonomous) function
$H_{(n)}$ is separable. In fact, referring to the results which can
be found in \cite{Crampin2002} for example, we can state more
precisely that we are looking at a separable system of St\"ackel
type and that the eigenfunctions $u^a(y,x)$ of $\jb$ are orthogonal
separation coordinates. It should then be true indeed that in those
coordinates we have $n$ time-independent quadratic integrals in
involution. Note that for the function $H_{(n)}$, for example, this
time-independence means among other things (see the expression
(\ref{newerH(n)})) that the function $W_{(n)} - (\tr \jb)W^1$, when
passing from the coordinates $(y,x)$ to the coordinates
$(y,u(y,x))$, should become independent of the $y$-variables.

To give some further backing for all these rather subtle properties,
we carry out more explicit calculations in the Appendix. The
programme there is that we complete the recursive scheme started in
Section~3 by proving, in analogy with the expressions (\ref{An},
\ref{An1}), explicit results for the different parts of all the
$A_{(i)}$ tensors and computing explicit expressions for the
corresponding functions $H_{(i)}$.

\section{Examples}

For examples which exhibit all features of our generalization, we
probably should think of systems with at least three degrees of
freedom and a non-Euclidean metric. But going into the details of
such applications would substantially add to the length of this
already rather long paper, so we will keep this for a separate study
and limit ourselves here to simpler two-dimensional situations. The
first example is inspired by an integrability study in \cite{Sar} on
what were called `generalized H\'enon-Heiles systems', which are
\sode s of the form
\begin{eqnarray*}
\ddot{q}_1 &=& -c_1 q_1 + b q_1^2 - aq_2^2,\\
\ddot{q}_2 &=& -c_2 q_2 - 2 m q_1q_2.
\end{eqnarray*}
Compared to older case studies of integrability of H\'enon-Heiles
systems, the generalization comes from the extra parameters $a$ and
$m$, which are motivated by allowing a Lagrangian description of the
system in which the Hessian of the Lagrangian need not be (a
constant multiple of) the unit matrix. The investigation carried out
in \cite{Sar} mainly consisted in looking for all possible parameter
cases for which the system has two independent quadratic integrals.
It led to the identification of three new cases, which are in some
sense degenerate cases, because either $a$ or $m$ is zero, meaning
that the Lagrangian one originally thought of is degenerate. As a
result, there is no corresponding Hamiltonian which in the standard
cases is always available as the first of two integrals in
involution. A subsidiary question then was: to what extent can the
two integrals in those degenerate cases be understood as being in
involution, and it was argued that this can be resolved by
constructing, in principle, a suitably adapted non-standard Poisson
structure. As the equations in the case that either $a$ or $m$ is
zero clearly exhibit partial decoupling, there is a good chance that
those degenerate cases actually fit within our present theory, which
is what we will discuss now.

Consider the case where $m=0$ and $b=0$, so that the system reduces
to
\begin{eqnarray*}
\ddot{q}_1 &\!=\!& -c_1 q_1 - aq_2^2,\\
\ddot{q}_2 &\!=\!& -c_2 q_2.
\end{eqnarray*}
These equations are of course easy to solve without further ado, but
they must serve here in the first place to illustrate various
aspects of our theory. The second equation plays the role of driving
equation. For consistency with the notations in the preceding
sections we rename the variables as $q_2=y$, $q_1=x$ and write
\begin{eqnarray*}
\ddot{y} &=& -c_2 y,\\
\ddot{x} &=& -c_1 x - ay^2.
\end{eqnarray*}
Putting $p_y = \dot{y}$ and $p_x=\dot{x}$ the two quadratic
integrals read
\begin{eqnarray*}
F_2 &=& \onehalf p_y^2 + \onehalf c_2 y^2\\
F_1 &=& \onehalf (c_1 - 4c_2) p_x^2 + 2ayp_xp_y - 2 a xp_y^2
+ \onehalf c_1 (c_1-4c_2)x^2 \\
&& \qquad \mbox{} + a(c_1-2c_2) x y^2 +\onehalf a^2 y^4.
\end{eqnarray*}
Obviously $F_2$ is a Hamiltonian for the driving equation and is
likely to be identifiable with the function $H_{(2)}$ in the theory
(see (\ref{H(n+1)}) knowing that $n=1$ here). The idea now is the
following. Since $F_1$ is a first integral of the complete system,
its quadratic part identifies a Killing tensor $A$. Looking at the
tensor $J$ of which $A$ is the cofactor, this may or may not be a
scKt in general. It was shown in \cite{Lundmark2} however that this
will always be the case for the Euclidean metric in dimension 2. If
the component $J_1$ is nonsingular therefore, we must be in a
situation covered by our present theory and all the features we
discussed should apply, with $F_1=H_{(1)}=H$, the Hamiltonian of the
quasi-Hamiltonian representation (\ref{quasiHam}). The extra
assumption that the driven system should have a genuine potential,
parametrically depending on the $y$-variables, is also automatically
satisfied here by dimension.

From the expression of $F_1$, we see that
\begin{equation}
A= \left(\begin{array}{cc}
-4ax & 2ay\\
2ay & c_1 - 4c_2
\end{array}\right).
\label{Aex1}
\end{equation}
and $A=\cof J$, where $J$ is the tensor with the following matrix
components
\begin{equation}
J= \left(\begin{array}{cc}
c_1-4c_2 & -2ay\\
-2ay & -4ax
\end{array}\right).
\label{Jex1}
\end{equation}
This tensor is indeed a scKt with respect to the Euclidean metric,
and its $J_1$ component is nonsingular (assuming $c_1 \neq 4c_2$),
so we are in business and $A^1 = \cof J_1 =1$. For completeness,
observe that
\[
\mu = - c_2y dy - (c_1x + ay^2)dx.
\]
in this example and that $A\mu = -dW$ indeed, with
\[
W=\onehalf c_1 (c_1-4c_2)x^2 + a(c_1-2c_2) x y^2 +\onehalf a^2 y^4.
\]
The function $h$ for the standard Hamiltonian representation of the
driven equation (with parameter $y$) is given by
\[
h= \onehalf p_x^2 + a y^2 x +  \onehalf c_1 x^2.
\]

Let us now further illustrate the subtleties of the theory, as
explained in Sections 4 to 6. We start by computing the function
$\tilde{h}$ as defined in (\ref{tildeh}). The tensor $\jb$ here
reads
\begin{equation}
\jb = - \left(4ax + \frac{4a^2y^2}{c_1-4c_2}\right)\fpd{}{x}\otimes
dx. \label{jbex1}
\end{equation}
and the new momenta defined in (\ref{tildep_i}, \ref{tildep_a})
become
\begin{eqnarray*}
p_y &=& (c_1-4c_2)\tilde{p}_y,\\
p_x &=& \tilde{p}_x -2ay \tilde{p}_y.
\end{eqnarray*}
In those new variables, the expression for $F_1$ becomes
\begin{eqnarray*}
F_1 &=& \onehalf (c_1 -4c_2)\tilde{p}_x^2 -
2a^2y^2(c_1-4c_2)\tilde{p}_y^2
-2ax(c_1-4c_2)^2 \tilde{p}_y^2 \\
&&\qquad + \onehalf c_1 (c_1-4c_2)x^2 + a(c_1-2c_2) x y^2 +\onehalf
a^2 y^4
\end{eqnarray*}
which is in agreement with (\ref{newH(n)}). For the computation of
$\tilde{h}$ on the other hand, we need the generating function
(\ref{F2}) of the time-dependent canonical transformation
$(x,p_x)\leftrightarrow (x,\tilde{p}_x)$. We see from (\ref{Jex1})
that $J^1_2= -2ay=\partial\psi/\partial x$ with $\psi=-2axy$. The
generating function $F$ thus reads
\[
F(x,\tilde{p}_x,t) = x\tilde{p}_x - 2 axy(t)\tilde{p}_y(t),
\]
and computing its partial time derivative involves making use of the
driving equation. The resulting expression for $\tilde{h}$ is found
to be
\begin{eqnarray*}
\tilde{h} =\onehalf \tilde{p}_x^2 + 2a^2y^2\tilde{p}_y^2 -
2ax(c_1-4c_2)\tilde{p}_y^2-2ay\tilde{p}_x\tilde{p}_y+a y^2x
+\onehalf c_1 x^2 + \frac{2ay^2x c_2}{c_1-4c_2}.
\end{eqnarray*}
One can verify that $\tilde{h}$ and $H_{(1)}=F_1$ indeed verify the
requirement (\ref{fulltildeh}) of Lemma~2 to within an additive
function of time, which is the function
\[
4a^2y^2\tilde{p}_y^2 -\frac{a^2y^4}{2(c_1-4c_2)}
\]
and of course can be ignored in writing down Hamilton's equations.
For the final canonical transformation to be applied to $\tilde{h}$,
we need the eigenfunction $u(y,x)$ of $\jb$, which is found to be
(see (\ref{jbex1}))
\[
u(y,x)=-\frac{4a^2y^2}{c_1-4c_2}-4ax.
\]
The time-dependent canonical transformation with generating function
$F(x,s,t) = s\,u(y(t),x)$ now transforms the relevant part of
$\tilde{h}$ into the function
\[
(\det J_1)^{-1}F_1 = 8a^2 s^2 + \onehalf u (c_1-4c_2)\tilde{p}_y^2 +
\onehalf u \frac{c_2}{c_1-4c_2} y^2 + \frac{c_1}{32a^2}u^2.
\]
It then easily follows by taking into account that $F_2$ is a
constant along solutions of the driving equation, that
\begin{equation*}
F_1 = 8 (c_1-4c_2)a^2 s^2 + u F_2
  +\frac{c_1(c_1-4c_2)}{32a^2}u^2,
\end{equation*}
which indeed no longer depends explicitly on time.

The other two degenerate cases in \cite{Sar} also fit within the
present theory to some extent, but are even more peculiar. One can
verify that the two integrals in those cases again can be understood
by the fact that the system is a driven cofactor system. This also
explains why the equations partially decouple. But it so happens
that the $J_1$ part of the scKt $J$ in those cases is zero, so that
we are not in the generic case of a nonsingular $J_1$, which was the
assumption in the preceding sections.

With a second simple example, we want to illustrate and test mainly
the beginning of our theory, before suitable coordinates for partial
decoupling have been identified, i.e.\ the situation covered by the
conditions of Definition~2. Consider the system
\begin{eqnarray*}
\ddot{q}_1 &=& 5q_1 - 4 q_2, \\
\ddot{q}_2 &=& q_1.
\end{eqnarray*}
The connection coefficients of the \sode\ connection are all zero
and the $(1,1)$ tensor $\Phi$ has the following matrix
representation
\[
\left( \Phi^i_j\right) = \left( -\fpd{f^i}{q^j}\right) = \left(
\begin{array}{cc} -5 & 4 \\-1 & 0 \end{array} \right).
\]
With such a simple, constant Jacobi endomorphism, finding a
distribution $K$ which entails submersiveness of the system, i.e.\
which satisfies the conditions (\ref{K}) of Definition~2, is simply
a matter of looking for a 1-dimensional eigenspace of $\Phi$. We
choose
\[
K = {\rm sp\,} \left\{ \fpd{}{q_1} + \fpd{}{q_2} \right\}.
\]
Note that submersiveness now is ensured and is not related to the
existence of a cofactor representation of our system. The problem of
detecting such a representation is quite interesting in its own
right. It bears to some extent resemblance to the inverse problem of
the calculus of variations, because there is a certain freedom in
selecting a multiplier matrix $g$ first. For example, making the
obvious choice of the unit matrix for $g$, or expressed in more
mechanical terms, associating the given equations with a the
standard kinetic energy $T=\onehalf(\dot{q}_1^2 + \dot{q}_2^2)$ will
not work! Indeed, one can verify that with this $g$, there is no
special conformal Killing tensor $J$ for which the right-hand sides
of the equations will satisfy the condition $D_J\mu=0$. On the other
hand, the following choice for $g$ turns out to be appropriate. Take
\[
\left(g_{ij}\right) = \left( \begin{array}{cc} 1 & -1 \\ -1&10
\end{array} \right),
\]
then, a scKt $J$ with respect to $g$ is easily found to be
\[
\left(J_{ij}\right) = \left( \begin{array}{cc} 2q_1 - 2q_2 & q_1+8q_2 \\
q_1+8q_2 & -4q_1 + 40 q_2 \end{array} \right).
\]
But more importantly, the nonconservative forces of our system are
now of the appropriate form. That is to say, multiplying the
right-hand sides of the equations with $g$, the 1-form $\mu$ is
found to be
\[
\mu = 4(q_1-q_2)dq_1 + (5q_1+4q_2)dq_2,
\]
and computing the cofactor tensor $A$ of $J$, one can verify that
$A\mu$ indeed is closed (or equivalently $D_J\mu=0$). We can now
compute the orthogonal complement of the distribution $K$ and see
whether the final conditions for a driven cofactor system are
verified. Since $g_{11}+g_{21}=0$, it easily follows that
\[
K^\perp = {\rm sp\,}\left\{ \fpd{}{q_1}\right\}.
\]
We have
\[
\DH{}\mu = 4dq_1\otimes dq_1 + 5 dq_1\otimes dq_2 - 4dq_2\otimes
dq_1 + 4 dq_2\otimes dq_2.
\]
It is clear then that $\DH{}\mu(K^\perp,K) \neq 0$, whereas
obviously $d\mu (K,K)=0$ by dimension. Hence, all requirements of
Definition~2 are met. Integrating the distributions $K$ and
$K^\perp$, suitable coordinates for the decoupling are found to be
$x=q_2,\ y=q_1 - q_2$ and the transformed system becomes
\begin{eqnarray*}
\ddot{y} &=& 4y, \\
\ddot{x} &=& x+y.
\end{eqnarray*}
We leave it to the reader to verify from here, exactly as we did
with the first example, that all the other features of our theory
hold true, in particular those related to the consecutive canonical
transformations relating the Hamiltonian $h$ of the driven system to
the quadratic integral $H_{(1)}$.

\section{Appendix}

Consider the recursive scheme for the tensors $A_{(i)}$, in
particular the relations (\ref{cfn}) for $i=1,\ldots,n$. We have
seen that taking $i=n$ led to determining relations for three of the
four parts of $A_{(n)}$ (see (\ref{An})), while the remaining block
$A_{(n)1}$ had to be obtained from taking subsequently $i=n-1$. This
procedure can be continued all the way down, and the essential
features of the recursion are captured in the following statement.

\begin{prop} Suppose that we know the tensor parts $A_{(i+1)21}$,
$A_{(i+1)12}$ and $A_{(i+1)2}$ and that they satisfy the identities
\begin{eqnarray}
&& J_1 A_{(i+1)12} + J_{12} A_{(i+1)2} \equiv 0 \equiv A_{(i+1)21}
J_1 + A_{(i+1)2} J_{21}, \label{id1-2} \\
&& J_{21} A_{(i+1)12} + J_2 A_{(i+1)2} \equiv A_{(i+1)21} J_{12} +
A_{(i+1)2} J_2. \label{id3}
\end{eqnarray}
Then the following are determining equations for the completion of
the construction of $A_{(i+1)}$ and for the next step in the
recursion
\begin{eqnarray}
A_{(i+1)1} &=& \Delta_{(i+1)} J_1^{-1} - J_1^{-1} J_{12}
A_{(i+1)21}, \label{next1} \\
A_{(i)21} &=& - J_{21} A_{(i+1)1} - J_2 A_{(i+1)21}, \label{next2} \\
A_{(i)12} &=& - A_{(i+1)1} J_{12} - A_{(i+1)12} J_2, \label{next3}
\\
A_{(i)2} &=& \Delta_{(i+1)} P_2 - J_{21} A_{(i+1)12} - J_2
A_{(i+1)2}. \label{next4}
\end{eqnarray}
Moreover, the newly obtained parts of $A_{(i)}$ will satisfy the
same three identities as those assumed for $A_{(i+1)}$, while the
completion of $A_{(i+1)}$ will give rise to the following
supplementary identity:
\begin{equation}
J_1 A_{(i+1)1} + J_{12} A_{(i+1)21} \equiv A_{(i+1)1} J_1 +
A_{(i+1)12} J_{21}. \label{id4}
\end{equation}
\end{prop}
\begin{proof}
Consider the recursion relation (\ref{cfn}) which in fact has two
parts. The idea is to compose each of those parts on the left and on
the right with one of the projectors $P_i$. This gives rise to a
total of 8 equations. For example, acting with $P_1$ on both sides
of the relation $JA_{(i+1)}+ P_2A_{(i)} = \Delta_{(i+1)}I_N$ (which
we can call a $(P_1,P_1)$ action for brevity) implies that we must
have
\[
J_1 A_{(i+1)1} + J_{12}A_{(i+1)21} = \Delta_{(i+1)} P_1,
\]
from which the determining equation (\ref{next1}) follows. Likewise,
a two-sided $(P_2,P_1)$ action will generate (\ref{next2}) and the
double $P_2$ action generates (\ref{next4}). The remaining
$(P_1,P_2)$ combination merely confirms the first of the assumed
identities (\ref{id1-2}). Starting from the other part in
(\ref{cfn}), the determining equation (\ref{next3}) will follow from
a $(P_1,P_2)$ action. The $(P_2,P_1)$ combination confirms the
second of the assumed identities in (\ref{id1-2}). The double $P_2$
action gives rise to another determining equation for $A_{(i)2}$,
which is consistent with the first one in view of the identity
(\ref{id3}). Finally the double $P_1$ action implies that, again
for consistency, $A_{(i+1)}$ must further satisfy the identity
(\ref{id4}). It is then a straightforward computation to verify that
the obtained blocks $A_{(i)21}$, $A_{(i)12}$ and $A_{(i)2}$ will
satisfy corresponding identities of the form (\ref{id1-2}) and
(\ref{id3}) in view of those assumed for $A_{(i+1)}$. This concludes
the full recursion step. \qed
\end{proof}

The identities which the different parts of all $A_{(i)}$ tensors
satisfy are important to get to a considerable simplification of the
recursive scheme. Indeed, as will be shown now, it turns out that
knowledge of the block $A_{(i)2}$ suffices to determine the three
other blocks of $A_{(i)}$. Moreover we can set up a recursive
scheme to determine $A_{(i)2}$ from $A_{(i+1)2}$ and we shall see
that this procedure brings the tensor $\jb$ back into the
spotlights.

It follows from the identities (\ref{id1-2}) that for each
$i=1,\ldots, n$:
\begin{eqnarray}
A_{(i)12} &=& - J_1^{-1} J_{12} A_{(i)2}, \label{Ai12} \\
A_{(i)21} &=& - A_{(i)2} J_{21}J_1^{-1}, \label{Ai21}
\end{eqnarray}
and subsequently from the defining relation (\ref{next1}) that
\begin{equation}
A_{(i)1} = \Delta_{(i)}J_1^{-1} + J_1^{-1}J_{12}A_{(i)2}J_{21}
J_1^{-1}. \label{Ai1}
\end{equation}
Making use of (\ref{Ai12}, \ref{Ai21}) in the identity (\ref{id3})
it is easy to see that this expresses the commutativity
\begin{equation}
\jb A_{(i)2} = A_{(i)2} \jb.  \label{jbA2}
\end{equation}
Finally, the determining equation (\ref{next4}) reduces to
\begin{equation}
A_{(i)2} = \Delta_{(i+1)} P_2 - \jb A_{(i+1)2}. \label{2recursion}
\end{equation}

\begin{lemma} The blocks $A_{(i)2} = P_2 \circ A_{(i)}\circ P_2$ of
the $A_{(i)}$ tensors which determine the quadratic part of the
first integrals $H_{(i)}$ are recursively given by
\begin{equation}
A_{(i)2} = \Delta_{(i+1)}P_2 + \sum_{j=1}^{n-i} (-1)^j
\Delta_{(j+i+1)}\jb^j, \quad i=1,\ldots, n-1, \label{A2}
\end{equation}
and all other parts of $A_{(i)}$ follow from $A_{(i)2}$.
\end{lemma}
\begin{proof}
We know that $A_{(n)2} = \Delta_{(n+1)}P_2$, with
$\Delta_{(n+1)}=\det J_1$. The recursive relation (\ref{A2}) then
easily follows from (\ref{2recursion}) by induction. The last part
of the statement has already been proved above. \qed
\end{proof}

Obviously, we are now in a position to venture computing a more
explicit expression for the functions $H_{(i)}$ and it turns out
that it is most appropriate to do this in terms of the momenta
$\tilde{p}$ again.

\begin{prop} The quadratic integrals $H_{(i)}$ (for $i=1,\ldots, n$) are given by
\begin{equation}
H_{(i)} = \onehalf A_{(i)}^{ab}\tilde{p}_a\tilde{p}_b + \onehalf
\Delta_{(i)} J_1^{kl} \tilde{p}_k\tilde{p}_l + W_{(i)}, \label{H(i)}
\end{equation}
where
\begin{equation}
A_{(i)}^{ab} = \Delta_{(i+1)}g^{ab} + \sum_{j=1}^{n-i} (-1)^j
\Delta_{(j+i+1)}(\jb^j)^{ab}. \label{A2^ab}
\end{equation}
\end{prop}
\begin{proof}
Recall that $H_{(n)}$ has already been computed (see (\ref{H(n)}))
and is indeed of the form (\ref{H(i)}). From (\ref{H(a)}), we
further have that
\[
H_{(i)} = \onehalf A_{(i)}^{ab} p_ap_b + A_{(i)}^{ak}p_ap_k +
\onehalf A_{(i)}^{kl} p_kp_l + W_{(i)}.
\]
Observe first that raising indices in (\ref{2recursion}) gives rise
to the formula
\[
A_{(i)}^{ab}= \Delta_{(i+1)} g^{ab} - \jb^b_d A_{(i+1)}^{da},
\]
and it subsequently follows from (\ref{Ai12}) and (\ref{Ai1}) that
\begin{eqnarray*}
A_{(i)}^{ak} &=& {J_1^{-1}}^k_l J^l_b(\jb^b_dA_{(i+1)}^{da}
- \Delta_{(i+1)} g^{ab}) \\
A_{(i)}^{kl} &=& \Delta_{(i)}{J_1^{-1}}^{kl} -
{J_1^{-1}}^l_jJ^j_b\jb^b_d {A_{(i+1)}}^d_cJ^c_m {J_1^{-1}}^{mk} +
\Delta_{(i+1)}{J_1^{-1}}^l_jJ^j_cJ^c_m{J_1^{-1}}^{mk}.
\end{eqnarray*}
We now compute each of the quadratic parts of $H_{(i)}$ in terms of
the $\tilde{p}$, using (\ref{tildep_i}, \ref{tildep_a}). The first
term, for example, becomes
\begin{eqnarray*}
\lefteqn{\onehalf A_{(i)}^{ab} p_ap_b = \onehalf A_{(i)}^{ab}
\tilde{p}_a\tilde{p}_b + \Delta_{(i+1)}J^{bj}\tilde{p}_b\tilde{p}_j
- \jb^b_d A_{(i+1)}^{da}J^j_b\tilde{p}_a\tilde{p}_j } \\[1mm]
& \hspace{2cm} + \onehalf
\Delta_{(i+1)}J^{bj}J^k_b\tilde{p}_j\tilde{p}_k - \onehalf \jb^b_d
A_{(i+1)}^{da} J^j_aJ^k_b \tilde{p}_j\tilde{p}_k .
\end{eqnarray*}
It so happens that only the first term on the right in this
expression survives, in other words all the other terms in the end
cancel out when we proceed in the same way with the two other
quadratic terms of $H_{(i)}$. The formula (\ref{H(i)}) then readily
follows, while (\ref{A2^ab}) merely is the contravariant form of
(\ref{A2}), and will turn out to be useful further on. \qed
\end{proof}

As was the case with $H_{(n)}$, we recognize in the expression for
the other $H_{(i)}$ part of the constant  $E^1$ of the driving
system. Explicitly, with the help of (\ref{E1}), it is clear that
along solutions of the driving system, the $H_{(i)}$ can be written
as
\begin{equation}
H_{(i)} = \onehalf A_{(i)}^{ab}\tilde{p}_a\tilde{p}_b +
\frac{\Delta_{(i)}}{\det J_1}(E^1-W^1) + W_{(i)}, \label{H(i)2}
\end{equation}
whereby we recall that $\det J_1=\Delta_{(n+1)}$. We shall now
finally illustrate that applying the parameter-dependent coordinate
change $(x^a,\tilde{p}_a)\leftrightarrow (u^a(x,y),s_a)$ turns the
expressions for the $H_{(i)}$ into functions which no longer depend
on the parameters $y^i$, and hence, in their interpretation of first
integrals of the driven system along solutions of the driving
system, become effectively time-independent quadratic integrals.

According to Proposition~3, the $u^a$ are eigenfunctions of $\jb$,
so that the coordinate expression of $\jb$, in the variables $(y,u)$
takes the simple form
\[
\jb = u^a \fpd{}{u^a}\otimes du^a.
\]
Going back to the result (\ref{Deltai}), it follows that
\[
\sum_{i=1}^{n+1} \Delta_{(i)}a^{i-1}= (\det J_1) \det(\jb + a I_n) =
(\det J_1) \prod_{b=1}^n(u^b+a),
\]
which in turn implies that
\begin{equation}
\Delta_{(i+1)} = (\det J_1)\, \sigma_{n-i}(u), \qquad i=0,\ldots, n,
\label{delta(i)}
\end{equation}
where $\sigma_{j}(u)$ denotes the elementary symmetric functions.

For information about the functional dependence of the other terms
in the potential part of the $H_{(i)}$ we go back to the double
cofactor representation of the overall system, which was the start
of the recursive scheme in Section~3. We know that with $A(a)$
representing the cofactor tensor of $J+aP_2$, the force terms $\mu$
of the overall system satisfy the relation $A(a)\mu=- dW(a)$, or
equivalently
\[
\det(J+aP_2) \mu = -(J+aP_2) dW(a).
\]
Projecting this relation under $P_1$ it follows that
\begin{equation}
\det(J+aP_2) \mu_1 = - J_1 P_1(dW(a)) - J_{21} P_2(dW(a)).
\label{aux}
\end{equation}
It is straightforward to compute that under the transformation
$(y,x) \leftrightarrow (y,u(x,y))$ the tensors involved acquire the
following coordinate expression
\[
P_1= \left(\fpd{}{y^i}+\fpd{u^b}{y^i}\fpd{}{u^b}\right) \otimes
dy^i, \qquad P_2= \fpd{}{u^a}\otimes du^a -
\fpd{u^a}{y^j}\fpd{}{u^a}\otimes dy^j,
\]
and
\[
J_1 = J^i_j \left(\fpd{}{y^i}+\fpd{u^b}{y^i}\fpd{}{u^b}\right)
\otimes dy^j, \qquad J_{21} = J^a_i\fpd{u^c}{x^a} \fpd{}{u^c}\otimes
dy^i.
\]
Concerning the left-hand side of (\ref{aux}) we have to remember
also that $(\det J_1) \mu_1 = - dW^1$, since the driving system is
of cofactor type with scKt $J_1$. It then readily follows, using the
property (\ref{auxprop4}), that the coordinate expression of
(\ref{aux}) reduces to
\[
\det(\jb+aP_2) \fpd{W^1}{y^k} = \fpd{W(a)}{y^k},
\]
and this for all $a$. The expansion (\ref{Deltai}) thus implies that
the functions $W_{(i)} - (\Delta_{(i)}/\det J_1) W^1$ do not depend
on the $y$ parameters, which together with (\ref{delta(i)}) confirms
our objective for the potential part in $H_{(i)}$.

It remains to look at the terms quadratic in the momenta. Expressed
in the momenta associated to the $u$-variables, they become
\[
\onehalf A^{ab}_{(i)}\fpd{u^c}{x^a}\fpd{u^d}{x^b}\,s_cs_d.
\]
We know that in particular $A^{ab}_{(n)}= (\det J_1)g^{ab}$ (cf.
(\ref{H(n)})). Since $\jb$ is a scKt with respect to $g_2$, it is
symmetric in its covariant or contravariant representation. But
$\jb$ is diagonal in the $u$ coordinates, hence the same is true for
the transformed $g_2$. Let us rely here further on the fact that we
proved by indirect means, mainly as a result of the simple Lemma~4,
that $H_{(n)}$ will be time-independent when expressed in the
$(u,s)$ variables. The net conclusion then is that the diagonal
elements of the transformed $g_2$ must be the product of $(\det
J_1)^{-1}(y)$ with a function depending on the $u$ variables only.
It subsequently follows from the explicit expression of
$A^{ab}_{(i)}$ in (\ref{A2^ab}) and the fact that the ratios
$\Delta_{(i)}/\det J_1$ also are functions of the $u$-variables only
(see (\ref{delta(i)})), that the quadratic part of all $H_{(i)}$
will indeed become time-independent as well, when these functions
are looked at as first integrals of the driven system along
solutions of the driving system.

To complete the picture it is perhaps worth repeating that the
final $H_{(1)}$ in that hierarchy is in fact the Hamiltonian $H$ in
the quasi-Hamiltonian representation of the full system and that
$(\det J_1)^{-1} A_{(1)2}$ is the cofactor of $\jb$. Acting with
$\jb$ on the expression (\ref{A2}) for $i=1$, we thus obtain a
polynomial expression satisfied by $\jb$. It is then comforting for
the internal consistency of our results that one can verify that
this identity indeed expresses the Cayley-Hamilton theorem applied
to $\jb$.

\subsubsection*{Acknowledgements} This work is part of the IRSES
project GEOMECH (nr. 246981) within the 7th European Community
Framework Programme. We are indebted to the referees for
constructive comments and suggestions.


{\footnotesize
\begin{thebibliography}{99}

\bibitem{Benenti}
S.\ Benenti, Inertia tensors and St\"ackel systems in the Euclidean
spaces, \emph{Rend.\ del Sem.\ Mat.\ Torino\/} \textbf{50} (1992)
315--341.

\bibitem{Benenti97}
S.\ Benenti, Intrinsic characterization of the variable separation
in the Hamilton-Jacobi equation, \emph{J.\ Math.\ Phys.\/} {\bf 38}
(1997) 6578--6602.

\bibitem{Benenti2}
S.\ Benenti, Special symmetric two-tensors, equivalent dynamical
systems, cofactor and bi-cofactor systems, \emph{Acta Appl.\
Math.\/} \textbf{87} (2005) 33--91.

\bibitem{Crampin2002}
M.\ Crampin, Conformal Killing tensors with vanishing torsion and
the separation of variables in the Hamilton-Jacobi equation,
\emph{Diff.\ Geom.\ Appl.\/} {\bf 18} (2002) 87--102.

\bibitem{CS2001}
M.\ Crampin and W.\ Sarlet, A class of nonconservative lagrangian
systems on Riemannian manifolds, \emph{J.\ Math.\ Phys.} \textbf{42}
(2001) 4313--4326.

\bibitem{CST}
M.\ Crampin, W.\ Sarlet and G.\ Thompson, Bi-differential calculi,
bi-Hamiltonian systems and conformal Killing tensors, {\em J.\
Phys.\ A:\ Math.\ Gen.\/} {\bf 33} (2000) 8755--8770.

\bibitem{FroNij}
A.\ Fr\"olicher and A.\ Nijenhuis, Theory of vector-valued
differential forms, {\em Proc.\ Ned.\ Acad.\ Wetensch. Ser.\ A\/}
{\bf 59} (1956) 338--359.

\bibitem{KT}
M.\ Kossowski and G.\ Thompson, Submersive second-order differential
equations, \emph{Math.\ Proc.\ Camb.\ Phil.\ Soc.\/} \textbf{110}
(1991) 207--224.

\bibitem{Lundmark2} H.\ Lundmark, Higher-dimensional integrable Newton systems with
quadratic integrals of motion, \emph{Studies in Appl.\ Math.\/}
\textbf{110} (2003) 257--296.

\bibitem{LW}
H.\ Lundmark and S.\ Rauch-Wojciechowski, Driven Newton equations
and separable time-dependent potentials, \emph{J.\ Math.\ Phys.}
\textbf{43} (2002) 6166--6194.

\bibitem{MarBlas}
K.\ Marciniak and M.\ B{\l}aszak, Non-Hamiltonian systems separable by
Hamilton-Jacobi method, \emph{J.\ Geom.\ Phys.} \textbf{58} (2008)
557–575.

\bibitem{Martinez}
E.\ Mart\'{\i}nez, Parallel transport and decoupling, In:
\emph{Applied differential geometry and mechanics}, Volume in honour
of the 60th birthday of Michael Crampin, W.\ Sarlet \& F.\ Cantrijn
(eds.) (Academia Press, Gent) (2003) 83--93.

\bibitem{MaCaSaI}
E.\ Mart\'{\i}nez, J.F.\ Cari\~{n}ena and W.\ Sarlet, Derivations of
differential forms along the tangent bundle projection, {\em Diff.\
Geometry and its Applications\/} {\bf 2} (1992) 17--43.

\bibitem{MaCaSaII}
E.\ Mart\'{\i}nez, J.F.\ Cari\~{n}ena and W.\ Sarlet, Derivations of
differential forms along the tangent bundle projection II, {\em
Diff.\ Geometry and its Applications\/} {\bf 3} (1993) 1--29.

\bibitem{MaCaSaIII}
E.\ Mart\'{\i}nez, J.F.\ Cari\~{n}ena and W.\ Sarlet, Geometric
characterization of separable second-order equations, {\it Math.\
Proc.\ Camb.\ Phil.\ Soc.\/} {\bf 113} (1993) 205--224.

\bibitem{WML}
S.\ Rauch-Wojciechowski, K.\ Marciniak and H.\ Lundmark,
Quasi-Lagrangian systems of Newton equations, \emph{J.\ Math.\
Phys.} \textbf{40} (1999) 6366--6398.

\bibitem{Sar}
W.\ Sarlet, New aspects of integrability of generalized
H\'enon-Heiles systems, \emph{J.\ Phys.\ A: Math.\ Theor.\/}
\textbf{24} (1991) 5245--5251.

\bibitem{SarBier}
W.\ Sarlet and W.\ Vanbiervliet, Geometric characterization of
driven cofactor systems, \emph{J.\ Phys.\ A: Math.\ Theor.\/}
\textbf{41} (2008) 042001 (10pp).



\end{thebibliography}

}
\end{document}